\documentclass[reqno,onecolumn,oneside]{paper}
\usepackage[english]{babel}
\usepackage{amsmath,amsthm,amsfonts,amssymb,latexsym,mathrsfs}
\usepackage[all]{xy}
\usepackage[sort]{cite}

\newtheorem{theorem}{Theorem}[section]
\newtheorem{proposition}[theorem]{Proposition}
\newtheorem{lemma}[theorem]{Lemma}
\newtheorem{corollary}[theorem]{Corollary}
\theoremstyle{definition}
\newtheorem{definition}[theorem]{Definition}
\newtheorem{exm}[theorem]{Example}

\theoremstyle{remark}
\newtheorem{remark}[theorem]{Remark}

\newcommand{\op}{\operatorname{op}}
\newcommand{\edm}{\operatorname{End}_\QQ}
\newcommand{\obj}{\operatorname{Obj}}

\newcommand{\id}{\operatorname{id}}

\newcommand{\frsl}{\operatorname{Free}_\SL}

\newcommand{\Quot}{\operatorname{Qt}}
\newcommand{\Cl}{\operatorname{Cl}}

\newcommand{\SL}{\mathcal{SL}}
\newcommand{\Q}{\mathcal{Q}}

\newcommand{\RL}{\mathcal{RL}}
\newcommand{\lMQ}{\mathcal{M}_\QQ^l}
\newcommand{\rMQ}{\mathcal{M}_\QQ^r}
\newcommand{\MQ}{\mathcal{M}_\QQ}
\newcommand{\cat}{\mathcal}
\newcommand{\Edm}{\mathbf{End}_\QQ}

\newcommand{\I}{\mathrm{I}}
\newcommand{\MM}{\mathbf{M}}
\renewcommand{\AA}{\mathbf{A}}
\newcommand{\NN}{\mathbf{N}}
\newcommand{\LL}{\mathbf{L}}

\newcommand{\QQ}{\mathbf{Q}}
\newcommand{\RR}{\mathbf{R}}

\renewcommand{\SS}{\mathbf{S}}
\newcommand{\Hom}{\mathbf{hom}}
\newcommand{\R}{\mathbb{R}}
\newcommand{\N}{\mathbb{N}}
\newcommand{\Z}{\mathbb{Z}}

\renewcommand{\wp}{\mathscr{P}}

\renewcommand{\phi}{\varphi}
\newcommand{\e}{\varepsilon}
\newcommand{\g}{\gamma}
\renewcommand{\d}{\delta}

\newcommand{\restr}{\upharpoonright}
\newcommand{\la}{\langle}
\newcommand{\ra}{\rangle}
\newcommand{\liff}{\Longleftrightarrow}
\newcommand{\upm}{_{\uparrow m}}

\newcommand{\under}{\backslash}
\newcommand{\ost}{{}_\star/}
\newcommand{\lto}{\longrightarrow}

\newcommand{\lmapsto}{\longmapsto}
\newcommand{\ust}{\under_\star}
\newcommand{\ov}{\overline}
\newcommand{\ul}{\underline}

\newcommand{\da}{\downarrow}
\newcommand{\ua}{\uparrow}

\newcommand{\0}{\mathbf{0}}

\def\amslatex\slash{{\protect\AmS-\protect\LaTeX}}

\begin{document} 

\title{Quantale Modules and their Operators, with Applications}

\author{Ciro Russo \\ \footnotesize{Dipartimento di Matematica ed Informatica -- Universit\`a di Salerno, Italy} \\ \footnotesize{Department of Mathematics -- Vanderbilt University, Nashville, USA} \\ \footnotesize{\texttt{cirusso@unisa.it}}}

\date{\today}
\maketitle

\begin{abstract}
The central topic of this work is the categories of modules over unital quantales. The main categorical properties are established and a special class of operators, called $\Q$-module transforms, is defined. Such operators~---~that turn out to be precisely the homomorphisms between free objects in those categories~---~find concrete applications in two different branches of image processing, namely fuzzy image compression and mathematical morphology.
\end{abstract}

\section{Introduction}

Actions of residuated structures on posets have been studied by many authors (see, for instance,~\cite{abramskyvickers, galatostsinakis, joyaltierney, mulvey, paseka, resende1, resende2, resende3, resende4, resende5, resendevickers, rosenthal, russo}), especially in connection with mathematical logic and topology. Such structures are involved in the development of several~---~rather recent~---~interesting theories including, for instance, foundational aspects of quantum mechanics (see, e.g.~\cite{mulvey}), linear logic~\cite{girard}, and abstract deductive systems for propositional logics~\cite{galatostsinakis}.

In particular, the structure of a quantale module (or, equivalently, of a complete poset subject to a biresiduated action from a complete residuated lattice) appears often in such works, although the basic properties of these categories still need to be investigated. A first aim of this paper is to fill this gap.

The second motivation that stimulated our investigation on quantale modules comes from the area of image processing. Indeed, in the literature of image processing, several suitable representations of digital images as $[0,1]$-valued maps are proposed. Such representations are the starting point for defining both compression and reconstruction algorithms based on fuzzy set theory~---~also called \emph{fuzzy algorithms}~---~and mathematical morphological operators, used for shape analysis in digital images.

Most of the fuzzy algorithms use a suitable pair of operators, one for compressing the image and the other one for approximating the original image starting from the compressed one; see, for instance,~\cite{dinolarusso, loiasessa, perfilieva}. The idea is similar to that of the so-called ``integral transforms'' in mathematical analysis: every map can be discretized by means of the \emph{direct transform} and then approximated through the application of a suitable \emph{inverse transform}. Usually a direct integral transform is defined as an integral-product composition; analogously, a fuzzy compression operator is defined as a join-product composition (where the product is actually a left-continuous triangular norm), and its inverse operator has the form of a meet-division composition, where the division is the residual operation of the same triangular norm.

In mathematical morphology, the operators of \emph{dilation} and \emph{erosion}~---~whose action on an image can easily be guessed by their names~---~that are translation invariant can be expressed, again, as join-product and meet-residuum compositions, respectively. All these methods can be placed under a common roof by essentially abstracting their shared properties. Indeed they are all examples of \emph{$\Q$-module transforms}, that we will define in Section~\ref{nucleisec} and that turn out to be precisely the homomorphisms between free $\Q$-modules.

The paper is meant to be as self-contained as possible and is organized as follows. The first two sections are preliminary and most of the results presented are known. Since quantale modules are basically sup-lattices subject to a biresiduated action from a quantale, in Section~\ref{resmapsec} we recall definitions and properties of residuated maps and sup-lattices, while Section~\ref{reslatquantsec} is dedicated to residuated lattices and quantales, the latter being essentially complete residuated lattices. It is important to underline immediately that, according to Theorem~\ref{quotmeetclosure}, residuated maps, sup-lattice morphisms and closure operators are pairwise in one-one correspondence. More precisely, sup-lattices and residuated maps between sup-lattices coincide; moreover, each closure operator gives rise to a sup-lattice morphism and each sup-lattice morphism, composed with its residual map, yield a closure operator.

In Section~\ref{mqbasnotsec} categorical properties of quantale modules are investigated. For the reader's convenience, it is divided into two subsections. The first contains basic definitions (objects, morphisms, subobjects, interval modules and so on) and results on quantale modules. In Subsection~\ref{freemqsub}, we characterize free objects, endow hom-sets with a sup-lattice structure and, lastly, show how products and coproducts are made.

We introduce $\Q$-module transforms and structural closure operators in Section~\ref{nucleisec}. Here we also show that $\Q$-module transforms are precisely morphisms of free modules~---~Theorem~\ref{repr}~---~and that structural closure operators are to $\Q$-module morphisms what closure operators are to sup-lattice morphisms, namely that each structural closure operator gives rise to a $\Q$-module morphism and vice versa. Morover, Theorem~\ref{repr} is extended in a suitable way to all $\Q$-modules (Theorem \ref{projrepr}).

%
Before concluding, in Section~\ref{appsec}, we show the relationships that link $\Q$-module transforms with fuzzy algorithms for image compression and reconstruction and with mathematical morphological operators. Also this section is divided into subsections. In Subsection~\ref{imagesec} the main idea underlying fuzzy algorithms for image compression is described. Subsection~\ref{matmorsec} contains the bases of mathematical morphology and a brief overview of its most important operators. Then, in Subsection~\ref{ordapprsec}, we discuss the unification of the aforementioned algorithms and operators in the framework of quantale modules.

\section{Residuated maps and sup-lattices}
\label{resmapsec}

Quantales and quantale modules are basically sup-lattices endowed with further operations whose behaviour can be described in the frameworks of residuation theory (see~\cite{blythjanowitz} for an introduction) and residuated lattices (for which the reader may refer to~\cite{jipsentsinakis}). Therefore in this section we recall some basic necessary definitions and results on residuated maps and sup-lattices.

\begin{definition}\label{residuated map}
Let $\la X, \leq \ra$ and $\la Y, \leq \ra$ be two posets. A map $f: X \lto Y$ is said to be \emph{residuated} provided there exists a map $g: Y \lto X$ such that, for all $x \in X$ and for all $y \in Y$, the following condition holds:
\begin{equation*}
f(x) \leq y \quad \iff \quad x \leq g(y).
\end{equation*}
It is immediate to verify that the map $g$ is uniquely determined; we call it the \emph{residual map} or the \emph{residuum} of $f$, and denote it by $f_*$. The pair $(f, f_*)$ is said to be \emph{adjoint}.
\end{definition}

Before discussing the basic properties of adjoint pairs, we recall that, if $\la X, \leq \ra$ is a poset, a map $\g: X \lto X$ is called a \emph{closure operator} if it is order preserving, extensive and idempotent, i.e. iff for all $x, y \in X$
\begin{enumerate}
\item[$(i)$]$x \leq y$ implies $\g(x) \leq \g(y)$,
\item[$(ii)$]$x \leq \g(x)$,
\item[$(iii)$]$\g \circ \g = \g$.
\end{enumerate}
Dually, a map $\d: X \lto X$ is called a \emph{coclosure operator}, or an \emph{interior operator}, if
\begin{enumerate}
\item[$(i)$]$\d$ is order preserving,
\item[$(ii)$]$\d(x) \leq x$ for all $x \in X$,
\item[$(iii)$]$\d \circ \d = \d$.
\end{enumerate}

The following result is a classical characterization of residuated maps:
\begin{theorem}\label{residchar}
Let $\la X, \leq \ra$ and $\la Y, \leq \ra$ be two posets, and $f: X \lto Y$. The following statements are equivalent:
\begin{enumerate}
\item[$(a)$]$f$ is residuated, with residual $f_*$;
\item[$(b)$]$f$ is monotone and for all $y \in Y$ there exists in $X$ the element
$$f_*(y) = \bigvee\{x \in X \mid f(x) \leq y\};$$
\item[$(c)$]$f$ is monotone and there exists a unique monotone map $f_*: Y \lto X$ such that
\begin{equation}\label{resid1}
f \circ f_* \leq \id_Y
\end{equation}
and
\begin{equation}\label{resid2}
f_* \circ f \geq \id_X;
\end{equation}
\item[$(d)$]$f$ is monotone and there exists a unique monotone map $f_*: Y \lto X$ such that
\begin{equation}\label{resid3}
f \circ f_* \circ f = f
\end{equation}
and
\begin{equation}\label{resid4}
f_* \circ f \circ f_* = f_*.
\end{equation}
\end{enumerate}
\end{theorem}

\begin{corollary}\label{resclosureinterior}
Let $\la X, \leq \ra$ and $\la Y, \leq \ra$ be two posets and let $(f, f_*)$ be an adjoint pair, with $f: X \lto Y$. Then
\begin{enumerate}
\item[$(i)$]$f_* \circ f$ is a closure operator over $X$;
\item[$(ii)$]$f \circ f_*$ is an interior operator in $Y$.
\end{enumerate}
\end{corollary}
\begin{proof}
It is enough to apply (\ref{resid1}--\ref{resid4}).
\end{proof}

Moreover we have:
\begin{proposition}\label{residprop}
Let $\la X, \leq \ra$ and $\la Y, \leq \ra$ be posets, and let $(f,f_*)$ be an adjoint pair, with $f: X \lto Y$. Then the following hold:
\begin{enumerate}
\item[$(i)$]$f$ preserves all existing joins, i.e. if $\{x_i\}_{i \in I}$ is a family of elements of $X$ such that there exists $\bigvee_{i \in I} x_i$, then also $\bigvee_{i \in I} f(x_i)$ exists and
$$f\left(\bigvee_{i \in I} x_i\right) = \bigvee_{i \in I} f(x_i);$$
\item[$(ii)$]$f_*$ preserves all existing meets, i.e. if $\{y_j\}_{j \in J}$ is a family of elements of $Y$ such that $\bigwedge_{j \in J} y_j$ exists, then also $\bigwedge_{j \in J} f_*(y_j)$ exists and
$$f_*\left(\bigwedge_{j \in J} y_j\right) = \bigwedge_{j \in J} f_*(y_j);$$
\item[$(iii)$]$f$ is surjective \quad $\iff$ \quad $f_*$ is injective \quad $\iff$ \quad $f \circ f_* = \id_Y$;
\item[$(iv)$]$f$ is injective \quad $\iff$ \quad $f_*$ is surjective \quad $\iff$ \quad $f_* \circ f = \id_X$.
\end{enumerate}
\end{proposition}

Let $\la X, \leq \ra$, $\la Y, \leq \ra$ and $\la W, \leq \ra$ be posets. A map $f: X \times Y \lto W$ of two variables is said to be \emph{biresiduated} if it is residuated with respect to each variable, i.e. if the following two conditions hold:
\begin{enumerate}
\item[-]for any fixed $\ov y$, there exists a map $g_{\ov y}: W \lto X$ such that
\begin{equation*}
f(x,\ov y) \leq w \qquad \iff \qquad x \leq g_{\ov y}(w);
\end{equation*}
\item[-]for any fixed $\ov x$, there exists a map $h_{\ov x}: W \lto Y$ such that
\begin{equation*}
f(\ov x,y) \leq w \qquad \iff \qquad y \leq h_{\ov x}(w).
\end{equation*}
\end{enumerate}

A poset $\la L, \leq \ra$ which admits arbitrary joins is called a \emph{sup-lattice}. A sup-lattice homomorphism is a map that preserves arbitrary joins. If $S \subseteq L$, the join over $S$ will be denoted indifferently by $\bigvee_{s \in S} s$ or $\bigvee S$. It is easily seen that a sup-lattice admits also arbitrary meets. Indeed, for any $S \subseteq L$, we can consider the set $S' = \{x \in L \mid s \leq x \ \forall s \in S\}$ of the lower bounds of $S$, and we have $\bigvee S' = \bigwedge S$. Then the category $\SL$ of sup-lattices is the one whose objects are complete lattices and morphisms are maps preserving arbitrary joins. If $\LL_1$ and $\LL_2$ are sup-lattices and $\phi: \LL_1 \lto \LL_2$ is a sup-lattice homomorphism, then clearly $\phi(\bot_1) = \phi\left({}^{\LL_1}\bigvee \varnothing\right) = {}^{\LL_2}\bigvee \varnothing = \bot_2$, where $\bot_i$ is the bottom element of $\LL_i$, $i = 1, 2$. Thus, for a sup-lattice $\LL$, we will use the notation $\LL = \la L, \vee, \bot \ra$.

For any sup-lattice $\LL = \la L, \vee, \bot \ra$ it is possible to define a dual sup-lattice in an obvious way: if we consider the opposite partial order $\geq$ (also denoted by $\leq^{\op}$), then $\LL^{\op} = \la L, \vee^{\op}, \bot^{\op} \ra = \la L, \wedge, \top \ra$~--- where $\top = \bigvee L$~--- is a sup-lattice and, clearly, $(\LL^{\op})^{\op} = \LL$.

It is also clear that, given two sup-lattices $\LL_1$ and $\LL_2$, $\hom_\SL(\LL_1,\LL_2)$ is a sup-lattice itself, with the order relation~--- and, therefore, the operations~--- defined pointwisely: $f \leq g \iff f(x) \leq g(x) \ \forall x \in L_1$. We will denote it by $\Hom_\SL(\LL_1,\LL_2)$.

If $\LL_1$ and $\LL_2$ are sup-lattices and $f \in \hom_\SL(\LL_1, \LL_2)$, $f$ is monotone and, for all $y \in L_2$, there exists $\bigvee\{x \in L_1 \mid f(x) \leq y\}$. Then, by Theorem~\ref{residchar}, $f$ is a residuated map whose residual, $f_*$, is defined by $f_*(y) = \bigvee\{x \in L_1 \mid f(x) \leq y\}$ for all $y \in L_2$. Moreover, $f_*$ is a homomorphism between the dual sup-lattices $\LL_2^{\op}$ and $\LL_1^{\op}$; to emphasize this fact, we will denote $f_*$ also by $f^{\op}$. For all $f, g \in \hom_\SL(\LL_1,\LL_2)$, we have:
\begin{enumerate}
\item[-] $(f^{\op})^{\op} = f$,
\item[-] $f \leq g \iff g_* \leq f_* \iff f^{\op} \leq^{\op} g^{\op}$,
\item[-] $(g \circ f)^{\op} = f^{\op} \circ g^{\op}$,
\end{enumerate}
hence
\begin{equation}\label{homsldual}
\Hom_\SL(\LL_1,\LL_2) \cong \Hom_\SL(\LL_2^{\op},\LL_1^{\op}).
\end{equation}
Conversely, if $f: \LL_1 \lto \LL_2$ is a residuated map between sup-lattices, it is a sup-lattice homomorphism between $\LL_1$ and $\LL_2$ by Proposition~\ref{residprop}$(i)$, and its residuum $f_*$ is a sup-lattice homomorphism from $\LL_2^{\op}$ to $\LL_1^{\op}$, by part $(ii)$ of the same proposition.

Then we have that sup-lattice homomorphisms coincide with residuated maps between sup-lattices. Furthermore, by Corollary~\ref{resclosureinterior}, if $f: \LL_1 \lto \LL_2$ is a residuated map, then $f_* \circ f: \LL_1 \lto \LL_1$ is a closure operator. Indeed it is also true that any closure operator $\g: \LL \lto \LL$ gives rise to a residuated map (i.e. a sup-lattice homomorphism) $f: x \in \LL \lmapsto \g(x) \in \g[\LL]$ (where the join ${}^\g\bigvee$ in $\g[\LL]$ is defined by ${}^\g\bigvee X = \g\left({}^\LL\bigvee X\right)$, for all $X \subseteq \g[\LL]$) whose residual is $f_*: y \in \g[\LL] \lmapsto \bigvee f^{-1}[\{y\}] \in \LL$, as we are going to see in Theorem~\ref{quotmeetclosure}.

Let $\Cl(\LL)$ be the set of all the closure operators over a given sup-lattice $\LL$; it is easily seen to be a poset with respect to the order relation defined pointwise: for all $\gamma, \d \in \Cl(\LL)$, $\gamma \leq \d \iff \gamma(x) \leq \d(x) \ \forall x \in L$. Furthermore, let us consider the set, $\Quot(\LL)$, of all the quotients of $\LL$ and the one, $\Lambda(\LL)$, of all the subsets of $\LL$ that are closed under arbitrary meets, both partially ordered by set inclusion. The following result proves that the posets $\Cl(\LL)^{\op}$, $\Lambda(\LL)$ and $\Quot(\LL)$ are isomorphic and, moreover, that any quotient of a sup-lattice is isomorphic~--- in a precise sense~--- to a subset of the same sup-lattice that is closed under arbitrary meets.

\begin{theorem}\label{quotmeetclosure}
For any sup-lattice $\LL$, $\Cl(\LL)^{\op}$, $\Quot(\LL)$ and $\Lambda(\LL)$ are isomorphic posets. Moreover, every element of $\Lambda(\LL)$ has a sup-lattice structure that is canonically isomorphic to a quotient of $\LL$.
\end{theorem}
\begin{proof}
For all $S \in \Lambda(\LL)$, we can define a map $\g_S: L \lto L$ by setting, for all $x \in L$, $\g_S(x) = \bigwedge\{y \in S \mid x \leq y\}$. It comes straightforwardly from its definition that $\g_S$ is monotone, extensive and idempotent; hence $\g_S \in \Cl(\LL)$ and $\g_S[L] = S$. Then we can define the map
$$\Gamma: \ S \in \Lambda(\LL) \ \lmapsto \ \g_S \in \Cl(\LL).$$
If $S, T \in \Lambda(\LL)$ are two different sets, then there exists $\ov x \in (S \setminus T) \cup (T \setminus S)$. If $\ov x \in S \setminus T$, then $\g_S(\ov x) = \ov x \neq \g_T(\ov x)$; analogously, if $\ov x \in T \setminus S$, then $\g_T(\ov x) = \ov x \neq \g_S(\ov x)$. Hence $S \neq T$ implies $\g_S \neq \g_T$, and $\Gamma$ is injective.

On the other hand, if $\g \in \Cl(\LL)$ and $S$ is an arbitrary subset of $\g[L]$, then $S$ is also a subset of $L$; thus there exists $z = \bigwedge S \in L$. Now, since $z \leq x$ for all $x \in S$, $\g(z) \leq \g(x) = x$ for all $x \in S$. Therefore $\g(z) \leq z$, whence $\g(z) = z \in S$. Then $\g[L]$ is closed under arbitrary meets. Moreover, if we assume that there exists $\ov x \in L$ such that $\g(\ov x) \neq \g_{\g[L]}(\ov x)$, then clearly $\g(\ov x) \gneqq \g_{\g[L]}(\ov x)$, by the definition of $\g_{\g[L]}$. So we have $\ov x \leq \g_{\g[L]}(\ov x)$ and $\g(\g_{\g[L]}(\ov x)) = \g_{\g[L]}(\ov x) < \g(\ov x)$, that contradicts the monotonicity of $\g$. It follows that $\g = \g_{\g[L]}$, for all $\g \in \Cl(\LL)$, and $\Gamma$ is surjective too, hence bijective. Now, if $S \subseteq T \in \Lambda(\LL)$, then clearly $\g_T \leq \g_S$ by the definition of such maps, so $\Gamma$ is order reversing. Analogously it is immediate to verify that also $\Gamma^{-1}$ is order reversing and therefore $\Gamma$ is an isomorphism between $\Lambda(\LL)$ and $\Cl(\LL)^{\op}$.

Now it is easy to see that, if $\g = \g_S$ is a closure operator, then $\mathbf{S} = \la S, \g \circ \vee, \g(\bot)\ra$ is a sup-lattice. Then, given a set $S \in \Lambda(\LL)$, we can define a map $\rho_S: L \lto S$ just by restricting the codomain of $\g_S$ to $S$:
\begin{equation}\label{reflection}
\rho_S: \ x \in L \ \lmapsto \ \g_S(x) = \bigwedge \{y \in S \mid x \leq y\} \in S.
\end{equation}
The map $\rho_S$, also called the \emph{reflection}\index{Sup-lattice!-- reflection} of $L$ over $S$, is the left inverse of the inclusion $\id_{L\restr S}: S \hookrightarrow L$ and is a surjective homomorphism from $\LL$ to $\mathbf{S}$; thus $\mathbf{S}$ is isomorphic to the quotient $\LL/\rho_S$ of $\LL$. In this way we define a map
$$\Pi: \ \mathbf{S} \in \Lambda(\LL) \ \lmapsto \ \LL/\rho_S \in \Quot(\LL),$$
and a map $\Theta = \Pi \circ \Gamma^{-1}$, from $\Cl(\LL)$ to $\Quot(\LL)$.

It is easily seen that $\Pi$ is an isomorphism and, consequently, $\Theta$ is an isomorphism between $\Cl(\LL)^{\op}$ and $\Quot(\LL)$; the theorem is proved.
\end{proof}

The following well-known result describes free objects in the category $\SL$.

\begin{proposition}\emph{\cite{joyaltierney}}\label{freesl}
For any set $X$, the free sup-lattice $\frsl(X)$ generated by $X$ is $\wp(\mathbf{X}) = \la \wp(X), \cup, \varnothing \ra$.
\end{proposition}

In particular, by Proposition~\ref{freesl}, the free sup-lattice over one generator $\wp 1$ is isomorphic to $\{\bot, \top\}$ and, clearly, $\wp 1^{\op} \cong \wp 1$. On the other hand, it is clear as well that, for an arbitrary sup-lattice $\LL$, a map $f: \wp1 \lto \LL$ is a sup-lattice homomorphism if and only if $f(\bot) = \bot_\LL$. Thus, using also (\ref{homsldual}) and the self-duality of $\wp1$, we have
$$\LL \cong \Hom_\SL(\wp1,\LL) \cong \Hom_\SL(\LL^{\op},\wp1)$$
and
$$\LL^{\op} \cong \Hom_\SL(\wp1,\LL)^{\op} \cong \Hom_\SL(\LL,\wp1),$$
hence
\begin{equation}\label{homp1l}
\LL \cong \Hom_\SL(\wp1,\LL) \cong \Hom_\SL(\LL^{\op},\wp1) \cong \Hom_\SL(\LL,\wp1)^{\op},
\end{equation}
\begin{equation}\label{homlp1}
\LL^{\op} \cong \Hom_\SL(\wp1,\LL)^{\op} \cong \Hom_\SL(\LL,\wp1) \cong \Hom_\SL(\LL^{\op},\wp1)^{\op}.
\end{equation}

\begin{proposition}\emph{\cite{joyaltierney}}\label{coprodsl}
Let $\{\LL_i\}_{i \in I}$ be a family of sup-lattices. Then the coproduct $\coprod_{i \in I} \LL_i$ is the product $\prod_{i \in I} \LL_i$ equipped with the maps $\mu_i: \LL_i \lto \prod_{i \in I} \LL_i$ that send each $x_i \in L_i$ into the family of $\prod_{i \in I} \LL_i$ in which all the elements are equal to $\bot$ except the $i$-th that is equal to $x_i$. Moreover, if $\pi_i: \prod_{i \in I} \LL_i \lto \LL_i$ is the canonical $i$-th projection, for all $i \in I$, $\pi_i \circ \mu_i = \id_{\LL_i}$.
\end{proposition}

\begin{proposition}\label{projsl}
Every free sup-lattice is projective. An object $\LL$ in $\SL$ is projective if and only if its dual $\LL^{\op}$ is injective.
\end{proposition}
\begin{proof}
Let $\wp(\mathbf{X})$ be a free sup-lattice, and let $\LL_1$ and $\LL_2$ be two sup-lattices such that there exists a surjective morphism $g: \LL_1 \lto \LL_2$ and a morphism $f: \wp(\mathbf{X}) \lto \LL_2$. Then $f$ is the unique homomorphism that extends the map $f_X: x \in X \lmapsto f(\{x\}) \in \LL_2$. If we consider the map $g_* \circ f_X: X \lto \LL_1$, we can extend it to a homomorphism $h: \wp(\mathbf{X}) \lto \LL_1$, and it is immediate to verify that $g \circ h = f$. Then any free sup-lattice is projective.

The second assertion is a trivial application of the self-duality of $\SL$.
\end{proof}

\section{Residuated lattices and quantales}
\label{reslatquantsec}

A binary operation $\cdot$ on a partially ordered set $\la P, \leq \ra$ is said to be \emph{residuated} iff there exist binary operations $\under$ and $/$ on $P$ such that for all $x, y, z \in P$,
\begin{equation*}
x \cdot y \leq z \quad \textrm{iff} \quad x \leq z/y \quad \textrm{iff} \quad y \leq x \under z.
\end{equation*}
The operations $\under$ and $/$ are referred to as the left and right \emph{residuals}, or \emph{divisions}, of $\cdot$, respectively. In other words, a residuated binary operation over $\la P, \leq \ra$ is a map from $P \times P$ to $P$ that is biresiduated. It follows from the results of Section~\ref{resmapsec} that the operation $\cdot$ is residuated if and only if it is order preserving in each argument and, for all $x, y, z \in P$, the inequality $x \cdot y \leq z$ has a largest solution for $x$ (namely $z/y$) and for $y$ (i.e. $x \under z$). In particular, the residuals are uniquely determined by $\cdot$ and $\leq$. The system $P = \la P, \cdot, \under, /, \leq \ra$ is called a \emph{residuated partially ordered groupoid} or \emph{residuated po-groupoid}.

In the situations where $\cdot$ is a monoid operation with a unit element $e$ and the partial order is a lattice order, we can add the monoid unit and the lattice operations to the similarity type to get an algebraic structure $\RR = \la R, \vee, \wedge, \cdot, \under, /, e \ra$ called a \emph{residuated lattice-ordered monoid} or \emph{residuated lattice} for short. It is not hard to see that $\RL$, the class of all residuated lattices, is a variety and the identities
\begin{equation*}
\begin{array}{lllll}
x \wedge (xy \vee z) / y \approx x, & \quad & x (y \vee z) \approx xy \vee xz, & \quad & (x / y)y \approx x \\
y \wedge x \under (xy \vee z) \approx y, & \quad & (y \vee z) x \approx yx \vee zx, & \quad & y(y \under x) \approx x, \\
\end{array}
\end{equation*}
together with monoid and lattice identities, form an equational basis for it.

In the category $\Q$ of quantales\footnote{In the literature, quantales are often defined as complete residuated po-groupoids, i.e. they are not unital. However, since we deal only with unital quantales, we will use this definition to avoid notational complications.}, $\obj(\Q)$ is the class of complete residuated lattices and the morphisms are the maps preserving products, the unit, arbitrary joins and the bottom element.

An alternative definition of quantale is the following

\begin{definition}\label{quantale}
A \emph{quantale} is an algebraic structure $\QQ = \la Q, \vee, \cdot, \bot, e \ra$ such that
\begin{enumerate}
\item[$(Q1)$]$\la Q, \vee, \bot \ra$ is a sup-lattice,
\item[$(Q2)$]$\la Q, \cdot, e \ra$ is a monoid,
\item[$(Q3)$]$x \cdot \bigvee\limits_{i \in I} y_i = \bigvee\limits_{i \in I} (x \cdot y_i)$ \ and \ $\left(\bigvee\limits_{i \in I} y_i\right) \cdot x = \bigvee\limits_{i \in I} (y_i \cdot x)$ \ for all $x \in Q$, $\{y_i\}_{i \in I} \subseteq Q$.
\end{enumerate}
$\QQ$ is said to be \emph{commutative} if so is the multiplication.
\end{definition}

The equivalence between the two definitions above is immediate to verify by means of the completeness of the lattice order in $\QQ$ and the properties of residuated maps. Indeed, since $\la Q, \vee, \bot \ra$ is a sup-lattice, it is possible to define the left and right residuals of $\cdot$:
\begin{equation*}
x \under y = \bigvee\{z \in Q \ \mid \ x \cdot z \leq y\},
\end{equation*}
\begin{equation*}
y / x = \bigvee\{z \in Q \ \mid \ z \cdot x \leq y\}.
\end{equation*}
Obviously, if $\QQ$ is commutative then the left and right divisions coincide:
\begin{equation*}
x \under y = y / x.
\end{equation*}

\begin{exm}\label{powermonoid}
Let $\AA = \la A, \cdot, e \ra$ be a monoid. We define, for all $X, Y \subseteq A$,
\begin{equation}\label{complexmultiplication}
\begin{array}{l}
X \cdot Y = \{x \cdot y \mid x \in X, y \in Y\},\\
X \cdot \varnothing = \varnothing \cdot X = \varnothing.
\end{array}
\end{equation}
It is immediate to verify that the structure $\wp(\AA) = \la \wp(A), \cup, \cdot, \varnothing, \{e\} \ra$ is a quantale, whose product is often called \emph{complex multiplication}.
\end{exm}

\begin{proposition}\label{freeqonm}
Let $\Q$ and $\cat M$ be, respectively, the categories of quantales and monoids. Then $\Q$ is a concrete category over $\cat M$. Moreover, given a monoid $\MM$, the free quantale over $\MM$ is $\wp(\MM)$.
\end{proposition}
\begin{proof}
The first part of the proposition is trivial. Indeed it is evident that the functor $U: \Q \lto \cat M$~--- defined as the functor that forgets the join and the bottom element~--- is faithful.

Now let us consider a monoid $\MM = \la M, \cdot, u \ra$, a quantale $\QQ = \la Q, \vee, \cdot, \bot, e \ra$ and a monoid morphism $f: \MM \lto U\QQ$. If we consider the singleton map $\sigma: x \in M \lto \{x\} \in \wp(M)$, it is obviously a monoid morphism and, if we set $h_f: X \in \wp(M) \lto \bigvee_{x \in X} f(x)$, the result follows as an easy application of Proposition~\ref{freesl}.
\end{proof}

Let $S$ be a non-empty set. It is well-known that the free monoid on $S$ is $\mathbf{S}^* = \la S^*, \cdot, \varnothing \ra$, where $S^* = \{x_1 \ldots x_n \mid n \in \N, x_i \in S\} \cup \{\varnothing\}$ and the product of two elements in $S^*$ is defined simply as their juxtaposition. With these notations, we have immediately the following result.

\begin{proposition}\label{freeqons}
Let $S$ be a non-empty set. The free quantale over $S$ is $\wp(\mathbf{S}^*)$.
\end{proposition}
\begin{proof}
If we consider the category $\cat M$ as a construct, with forgetful functor $U$, and $\Q$ as a concrete category over $\cat M$, with forgetful functor $U'$, it is immediate to verify that $U \circ U'$ is the underlying functor that makes $\Q$ a construct. Then $UU' \QQ = Q$, for any quantale $\QQ$. Now let $\QQ \in \Q$ and $f: S \lto Q$ be a map, and consider the following diagram
$$\xymatrix{
S \ar[rr]^{\id_{S^* \restr S}} \ar[ddrr]_f && S^* \ar[dd]^{Uh_f} \ar[rr]^\sigma && \wp(S^*) \ar[dd]^{UU'h_f'} \\
&&&&\\
	&& Q \ar[rr]_{\id_Q} && Q
},$$
where $\id_{S^* \restr S}$ is the inclusion map, $\sigma$ is the singleton map, $h_f$ is the unique monoid homomorphism that extends $f$ and $h_f'$ is the unique quantale homomorphism that extends $h_f$.

Then the diagram above is easily seen to be commutative, hence $h_f'$ is a homomorphism of quantales that extends $f$; the uniqueness of $h_f'$ comes from that of $h_f$. So the assertion is proved.
\end{proof}

Observe that, in $\wp(\mathbf{S}^*)$, $\varnothing$ is the bottom element and $\{\varnothing\}$ is the unit.

\begin{proposition}\label{aritq}
Let $\QQ$ be a quantale. Then, for any $x, y, z \in Q$ and $\{y_i\}_{i \in I} \subseteq Q$,
\begin{enumerate}
\item[$(i)$] $x \bot = \bot x = \bot$,
\item[$(ii)$] if $x \leq y$ then $xz \leq yz$ and $zx \leq zy$,
\item[$(iii)$] if $x \leq y$ then $x/z \leq y/z$, $z \under x \leq z \under y$, $z/y \leq z/x$ and $y \under z \leq x \under z$,
\item[$(iv)$] $(y/x)x \leq y$ and $x (x \under y) \leq y$,
\item[$(v)$] $x/e = e \under x = x$,
\item[$(vi)$] $x/\left(\bigvee_{i \in I} y_i\right) = \bigwedge_{i \in I}(x/y_i)$ and $\left(\bigvee_{i \in I} y_i\right) \under x = \bigwedge_{i \in I} (y_i \under x)$,
\item[$(vii)$] $\left(\bigwedge_{i \in I} y_i\right)/x = \bigwedge_{i \in I}(y_i/x)$ and $x \under \left(\bigwedge_{i \in I} y_i\right) = \bigwedge_{i \in I}(x \under y_i)$,
\item[$(viii)$] $y \under (x \under z) = xy \under z$ and $(z/y)/x = z/xy$.
\end{enumerate}
\end{proposition}
\begin{proof}
See Proposition 2.12 of~\cite{leethesis}.
\end{proof}

\section{The categories of Quantale Modules}
\label{mqbasnotsec}

\subsection{Basic definitions and properties}
\label{bdapsub}

\begin{definition}\label{modules}
Let $\QQ$ be a quantale. A (left) \emph{$\QQ$-module} $\MM$, or a \emph{module over $\QQ$}, is a sup-lattice $\la M, \vee, \bot \ra$ with an external binary operation, called \emph{scalar multiplication},
$$\star: (q,m) \in Q \times M \lmapsto q \star m \in M,$$
such that the following conditions hold:
\begin{enumerate}
\item[$(M1)$]$(q_1 \cdot q_2) \star m = q_1 \star (q_2 \star m)$, for all $q_1, q_2 \in Q$ and $m \in M$;
\item[$(M2)$]the external product is distributive with respect to arbitrary joins in both coordinates, i.e.
\begin{enumerate}
\item[$(i)$]for all $q \in Q$ and $\{m_i\}_{i \in I} \subseteq M$,
$$q \star {}^\MM\bigvee_{i \in I} m_i = {}^\MM\bigvee_{i \in I} q \star m_i,$$
\item[$(ii)$]for all $\{q_i\}_{i \in I} \subseteq Q$ and $m \in M$,
$$\left({}^\QQ\bigvee_{i \in I} q_i\right) \star m = {}^\MM\bigvee_{i \in I} q_i \star m,$$\
\end{enumerate}
\item[$(M3)$]$e \star m = m$.
\end{enumerate}
\end{definition}

Condition $(M2)$ can be expressed, equivalently, as follows:
\begin{enumerate}
\item[$(M2')$]The scalar multiplication is residuated with respect to the lattice order in $M$, i.e.
\begin{enumerate}
\item[$(i)$]for all $q \in Q$, the map
$$q^\star: m \in M \lmapsto q \star m \in M$$
is residuated,
\item[$(ii)$]for all $m \in M$, the map
$$^{\star}m: q \in Q \lmapsto q \star m \in M$$
is residuated.
\end{enumerate}
\end{enumerate}
Then, from $(M2')$ it follows that, for all $q \in Q$, there exists the residual map $(q^\star)_*$ of $q^\star$, and for all $m \in M$ there exists the residual map $(^\star m)_*$ of $^\star m$. Consequently Theorem~\ref{residchar} implies
\begin{equation}\label{qbulletsharp}
(q^\star)_*: m \in M \lmapsto \bigvee\{n \in M \mid q \star n \leq m\} \in M
\end{equation}
and
\begin{equation}\label{bulletmsharp}
(^\star m)_*: n \in M \lmapsto \bigvee\{q \in Q \mid q \star m \leq n\} \in Q.
\end{equation}
Condition (\ref{qbulletsharp}) defines another external operation over $M$:
$$\ust: (q,m) \in Q \times M \lmapsto q \ust m = (q^\star)_*(m) \in M.$$
Analogously, condition (\ref{bulletmsharp}) defines a map from $M \times M$ to $Q$:
$$\ost: (m,n) \in M \times M \lmapsto m \ost n = (^\star m)_*(n) \in Q.$$

The proof of the following proposition is straightforward from the definitions of $\star$, $\ust$ and $\ost$ and from the properties of quantales.

\begin{proposition}\emph{\cite{galatostsinakis}}\label{basicmq}
For any quantale $\QQ$ and any $\QQ$-module $\MM$, the following hold.
\begin{enumerate}
\item[$(i)$] The operation $\star$ is order-preserving in both coordinates.
\item[$(ii)$] The operations $\ust$ and $\ost$ preserve meets in the numerator; moreover, they convert joins in the denominator into meets. In particular, they are both order-preserving in the numerator and order reversing in the denominator.
\item[$(iii)$] $(m \ost n) \star n \leq m$.
\item[$(iv)$] $q \star (q \ust m) \leq m$.
\item[$(v)$] $m \leq q \ust (q \star m)$.
\item[$(vi)$] $(q \ust m) \ost n = q \under (m \ost n)$.
\item[$(vii)$] $((m \ost n) \star n) \ost n = m \ost n$.
\item[$(viii)$] $e \leq m \ost m$.
\item[$(ix)$] $(m \ost m) \star m = m$.
\end{enumerate}
\end{proposition}

\begin{exm}\label{freemodule}
Let $\QQ$ be a quantale and $X$ be an arbitrary non-empty set. We can consider the sup-lattice $\la Q^X, \vee^X, \bot^X \ra$, where $\bot^X$ is the $\bot$-constant function from $X$ to $Q$ and
$$\left(\bigvee_{i \in I}^X f_i\right)(x) = \bigvee_{i \in I} f_i(x) \quad \textrm{ for all } x \in X.$$
Then we can define a scalar multiplication in $Q^X$ as follows:
$$\star : (q,f) \in Q \times Q^X \lmapsto q \star f \in Q^X,$$
with the map $q \star f$ defined as $(q \star f)(x) = q \cdot f(x)$ for all $x \in X$.

It is clear that $Q^X$ is a left $\QQ$-module~--- denoted by $\QQ^X$~--- and, for all $q \in Q$, $f \in Q^X$ and $x \in X$, the following holds:
$$(q \ust f)(x) = q \under f(x).$$
\end{exm}

\begin{exm}\label{monoidaction}
Let $S$ be a set and $\AA = \la A, \cdot, e \ra$ be a monoid from which an action $\star$ on $S$ is defined. Then the sup-lattice $\wp(\SS) = \la \wp(S), \cup, \varnothing \ra$ is a module over the free quantale $\wp(\AA) = \la \wp(A), \cup, \cdot, \varnothing, \{e\} \ra$. Indeed it is easy to see that the module action $\star'$ of the monoid $\la \wp(A), \cdot, \{e\} \ra$ on the set $\wp(S)$, preserves arbitrary unions in both arguments.
\end{exm}

The definition and properties of right $\QQ$-modules are completely analogous. If $\QQ$ is commutative, the concepts of right and left $\QQ$-modules coincide and we will say simply $\QQ$-modules. If a sup-lattice $\MM$ is both a left $\QQ$-module and a right $\QQ'$-module~--- over two given quantales $\QQ$ and $\QQ'$~--- we will say that $\MM$ is a \emph{$(\QQ,\QQ')$-bimodule} if the following associative law holds:
\begin{equation}\label{bimodule}
(q \star_l m) \star_r q' = q \star_l (m \star_r q'), \quad \textrm{for all } \ m \in M, \ q \in Q, \ q' \in Q',
\end{equation}
where $\star_l$ and $\star_r$ are~--- respectively~--- the left and right scalar multiplications.

It is easy to prove, by means of Proposition~\ref{residprop}, that if $\QQ$ is a quantale and $\MM$ is a left (respectively: right) $\QQ$-module, the dual sup-lattice $\MM^{\op}$ is a right (resp.: a left) $\QQ$-module, with the external multiplication $\ust$ defined by condition (\ref{qbulletsharp}).

\begin{definition}\label{qmhomo}
Let $\QQ$ be a quantale and $\MM_1, \MM_2$ be two left $\QQ$-modules. A map $f: M_1 \lto M_2$ is a $\QQ$-module homomorphism if $f\left({}^{\MM_1}\bigvee_{i \in I} m_i\right) = {}^{\MM_2}\bigvee_{i \in I} f(m_i)$ for any family $\{m_i\}_{i \in I} \subseteq \MM_1$, and $f(q \star_1 m) = q \star_2 f(m)$, for all $q \in Q$ and $m \in M_1$, where $\star_i$ is the external product of $\MM_i$, for $i = 1,2$. The definition of right $\QQ$-module homomorphism is analogous
\end{definition}

Thus, given a quantale $\QQ$, the categories $\lMQ$ and $\rMQ$ have, respectively, left and right $\QQ$-modules as objects, and left and right $\QQ$-module homomorphisms as morphisms. If $\QQ$ is commutative, $\lMQ$ and $\rMQ$ coincide, and we denote such a category by $\MQ$.
\begin{remark}\label{notation}
Henceforth, in all the definitions and results that can be stated both for left and right modules, we will refer generically to ``modules''~--- without specifying left or right~--- and we will use the notations of left modules.
\end{remark}

\begin{proposition}\label{dualhomomq}
Let $\QQ$ be a quantale, $\MM_1$, $\MM_2$ be two $\QQ$-modules and $f: \MM_1 \lto \MM_2$ be a homomorphism. Then $f$ is a residuated map and the residual map $f_*: M_2 \lto M_1$ is a $\QQ$-module homomorphism between $\MM_2^{\op}$ and $\MM_1^{\op}$.
\end{proposition}
\begin{proof}
Since $f$ is a $\QQ$-module homomorphism, a fortiori it is a sup-lattice homomorphism, hence a residuated map and, as we proved in Section~\ref{resmapsec}, $f_*$ is a sup-lattice homomorphism from $\MM_2^{\op}$ to $\MM_1^{\op}$. What we need to prove, now, is that $f_*\left(q \under_{\star_2} m_2\right) = q \under_{\star_1} f_*(m_2)$. Let $m \in M_1$; we have
$$\begin{array}{llllll}
& m \leq f_*\left(q \under_{\star_2} m_2\right) & \ \liff \ & f(m) \leq q \under_{\star_2} m_2 & \ \liff \ & q \star_2 f(m) \leq m_2 \\
\liff \ & f(q \star_1 m) \leq m_2 & \ \liff \ & q \star_1 m \leq f_*(m_2) & \ \liff \ & m \leq q \under_{\star_1} f_*(m_2), \\
\end{array}$$
and the result follows.
\end{proof}

If $\MM_1$ and $\MM_2$ are left modules, then their duals are right modules and vice versa. Then it follows from Proposition~\ref{dualhomomq} that~--- as well as for sup-lattices~--- there is a one-one correspondence between the hom-set of two given modules and the hom-set of their dual modules. We will return to these properties (and clarifying them) in the next subsection.

Let $\QQ$ be a quantale and let $\MM = \la M, \vee, \bot \ra$ be a $\QQ$-module. A \emph{$\QQ$-submodule} $\NN$ of $\MM$ is a sup-sublattice $\la N, \vee, \bot \ra$ of $\la M, \vee, \bot \ra$ that is stable with respect to the external product of $\MM$. It is easy to verify that, for any family $\{\NN_i\}_{i \in I}$ of $\QQ$-submodules of $\MM$, $\left\la \bigcap_{i \in I} N_i, \vee, \bot \right\ra$ is still a $\QQ$-submodule of $\MM$. Thus, given an arbitrary subset $S$ of $M$, we define the $\QQ$-submodule $\left\la \la S \ra, \vee, \bot \right\ra$ \emph{generated by $S$} as the intersection of all the $\QQ$-submodules of $\MM$ containing $S$. Vice versa, given a submodule $\NN$ of $\MM$, we will say that a subset $S$ of $M$ is a \emph{system of generators} for $\NN$~--- or that $S$ \emph{generates} $\NN$~--- if $\NN = \la S \ra$.

If $\{\MM_i\}$ is a family of $\QQ$-modules, $\MM$ is a $\QQ$-module and $X$ is a non-empty set, the \emph{product} $\prod_{i \in I} \MM_i = \la \prod_{i \in I} M_i, \vee, (\bot_i)_{i \in I}) \ra$ of the family $\{\MM_i\}_{i \in I}$, and $\MM^X = \la M^X, \vee^X, \bot^X \ra$ are clearly $\QQ$-modules with the operations defined pointwise. $\MM^X$ is also called the \emph{power module} of $\MM$ by $X$.

\begin{proposition}\label{chsubgen}
Let $\QQ$ be a quantale, $\MM$ a $\QQ$-module and $S \subseteq M$. Then $\la S \ra = \bigvee Q \star S$, where
$$\bigvee Q \star S = \left\{\bigvee_{x \in S} q_x \star x \ \Big| \ \{q_x\}_{x \in S} \in Q^S\right\}.$$
\end{proposition}
\begin{proof}
First we observe that $\bot = \bigvee_{x \in S} \bot \star x \in \bigvee Q \star S$; now let $\{y_i\}_{i \in I}$ be an arbitrary family of elements of $\bigvee Q \star S$, with set of indices $I$. Then, for all $i \in I$, $y_i = \bigvee_{x \in S} q_x^i \star x$, for a suitable family $\{q_x^i\}_{x \in S}$ of elements of $Q$. We have
$$\bigvee_{i \in I} y_i = \bigvee_{i \in I} \left(\bigvee_{x \in S} q_x^i \star x\right) = \bigvee_{x \in S} \left(\bigvee_{i \in I} q_x^i \star x\right) = \bigvee_{x \in S} \left(\bigvee_{i \in I} q_x^i\right) \star x \in \bigvee Q \star S.$$
Moreover, for all $y = \bigvee_{x \in S} q_x \star x \in \bigvee Q \star S$ and for all $q \in Q$,
$$q \star y = q \star \bigvee_{x \in S} q_x \star x = \bigvee_{x \in S} q \star (q_x \star x) = \bigvee_{x \in S} (q \cdot q_x) \star x \in \bigvee Q \star S;$$
hence $\bigvee Q \star S$ is a $\QQ$-submodule of $\MM$ and it is clear that $S \subseteq \bigvee Q \star S$. Therefore $\la S \ra \subseteq \bigvee Q \star S$.

On the other hand, if $\NN$ is a $\QQ$-submodule of $\MM$ containing $S$, if we fix an element $q_x \in Q$ for each $x \in S$, the scalar product $q_x \star x$ must be in $N$ for all $x \in S$; thus also $\bigvee_{x \in S} q_x \star x \in N$; hence $\bigvee Q \star S \subseteq \NN$. The arbitrary choice of $\NN$, among the submodules of $\MM$ containing $S$, ensures that $\bigvee Q \star S \subseteq \la S \ra$; the result follows.
\end{proof}

Given a quantale $\QQ$, a $\QQ$-module is called \emph{cyclic} if it is generated by a single element $v$. According to the notation introduced in Proposition~\ref{chsubgen}, a cyclic module generated by a certain $v$, will be also denoted by $\QQ \star v$. 

\begin{lemma}\emph{\cite{galatostsinakis}}\label{cyclicchar}
A $\QQ$-module $\MM$ is cyclic with generator $v$ iff $(m \ost v) \star v = m$, for all $m \in M$.
\end{lemma}

\begin{definition}\label{intqmodule}
Let $\QQ$ be a quantale, and let $\MM = \la M, \vee, \bot \ra$ be a $\QQ$-module and $m$ be a fixed element of $M$. If we consider the set $M\upm = \{n \in M \mid m \leq n\}$, we can endow such a set with the structure of a $\QQ$-module. Indeed it is clear that $M\upm$ is closed both under arbitrary joins and meets; on the other hand, its bottom element is $m$ and we define the external operation $\star\upm$ as
\begin{equation*}
q \star\upm n = m \vee q \star n, \qquad \textrm{for all } n \in M\upm.
\end{equation*}
It is easy to verify that $\MM\upm = \la M\upm, \vee, m \ra$ is a $\QQ$-module, with such an external operation. We will call it the (upper) \emph{interval $\QQ$-module} determined by $m$ in $\MM$. For $M\upm$, we will use also the notation $[m, \top]$, indifferently.
\end{definition}

\subsection{Free modules, hom-sets, products and coproducts}
\label{freemqsub}

In the present subsection, we investigate several basic constructions and properties of the categories of quantale modules. According to Remark~\ref{notation}, in all the definitions and statements regarding modules over a non-commutative quantale $\QQ$, whenever we say simply $\QQ$-module or write $\MQ$, we mean that the definition or the result holds for both left and right modules (suitably reformulated, where necessary).

\begin{proposition}\label{freemq}
For any set $X$, the free $\QQ$-module generated by $X$ is the function module $\QQ^X = \la Q^X, \vee^X, \bot^X \ra$, with join and scalar multiplication defined pointwise, equipped with the map $\chi: x \in X \lmapsto \chi_x \in Q^X$, where $\chi_x$ is defined, for all $x \in X$, by
\begin{equation}\label{chi}
\chi_x(y) = \left\{\begin{array}{ll} \bot & \textrm{if } y \neq x \\ e & \textrm{if } y = x\end{array}\right..
\end{equation}
\end{proposition}
\begin{proof}
Let $\MM = \la M, \vee, \bot \ra$ be any $\QQ$-module and $f: X \lto M$ be an arbitrary map. We shall prove that there exists a unique $\QQ$-module morphism $h_f: \QQ^X \lto \MM$ such that $h_f \circ \chi = f$. First observe that, for any $\alpha \in Q^X$, $\alpha = \bigvee_{x \in X} \alpha(x) \star \chi_x$; then let us set $h_f(\alpha) = \bigvee_{x \in X} \alpha(x) \star f(x)$. For any fixed $\ov x \in X$, $(h_f \circ \chi)(\ov x) = h_f(\chi_{\ov x}) = \bigvee_{x \in X} \chi_{\ov x}(x) \star f(x) = e \cdot f(\ov x) \vee \bigvee_{x \in X \setminus \{\ov x\}} \bot \cdot f(x) = f(\ov x)$, hence $h_f \circ \chi = f$.

Let now $(\alpha_i)_{i \in I}$ be a family of elements of $Q^X$ and observe that $\bigvee_{i \in I} \alpha_i = \bigvee_{i \in I} \left(\bigvee_{x \in X} \alpha_i(x) \star \chi_x\right) = \bigvee_{x \in X} \left(\bigvee_{i \in I} \alpha_i(x) \star \chi_x \right)$. Then
$$\begin{array}{lll}
h_f\left(\bigvee_{i \in I} \alpha_i\right) &=& h_f \left(\bigvee_{i \in I} \left(\bigvee_{x \in X} \alpha_i(x) \star \chi_x \right) \right) \\
 &=& h_f\left(\bigvee_{x \in X} \left(\bigvee_{i \in I} \alpha_i(x) \star \chi_x \right) \right) \\
 &=& h_f\left(\bigvee_{x \in X} \left(\bigvee_{i \in I} \alpha_i(x) \right) \star \chi_x \right) \\
 &=& \bigvee_{x \in X} \left(\bigvee_{i \in I} \alpha_i(x)\right) \star f(x) \\
 &=& \bigvee_{x \in X} \left(\bigvee_{i \in I} \alpha_i(x) \star f(x)\right) \\
 &=& \bigvee_{i \in I} \left(\bigvee_{x \in X} \alpha_i(x) \star f(x)\right) \\
 &=& \bigvee_{i \in I} h_f(\alpha_i), \\
\end{array}$$
thus $h_f$ preserves arbitrary joins; the proof of the fact that it also preserves scalar multiplication is straightforward, and therefore $h_f$ is a $\QQ$-module homomorphism. Moreover, the uniqueness of $h_f$ can be proved exactly as in the proof of Proposition~\ref{freesl}.
\end{proof}
Obviously, every $\Q$-module is homomorphic image of a free module.

\begin{definition}\label{hommq}
Given $\Q$-modules $\MM$ and $\NN$, we define, on $\hom_\QQ(\MM,\NN)$, the following operations and constants:
\begin{enumerate}
\item[-]for all $\{h_i\}_{i \in I} \subseteq \hom_\QQ(\MM,\NN)$, the homomorphism $\bigsqcup_{i \in I} h_i$ is defined by $\left(\bigsqcup_{i \in I} h_i\right)(x) = \bigvee_{i \in I} h_i(x)$, for all $x \in M$,
\item[-]let $\bot^\bot$ and $\top^\top$ be the maps defined, respectively, by $\bot^\bot(x) = \bot_N$ and $\top^\top(x) = \top_N$, for all $x \in M$
\end{enumerate}
and, if $\QQ$ is commutative,
\begin{enumerate}
\item[-]for all $q \in Q$ and $h \in \hom_\QQ(\MM,\NN)$, let $q \diamond h$ be the map defined by $(q \diamond h)(x) = q \star h(x) = h(q \star x)$, for all $x \in M$.
\end{enumerate}
It is easy to see that $\Hom_\QQ(\MM,\NN) = \la \hom_\QQ(\MM,\NN), \sqcup, \bot^\bot \ra$ is a sup-lattice and, if $\QQ$ is a commutative quantale, it is a $\QQ$-module with the external multiplication $\diamond$. If $\NN = \MM$, the sup-lattice (or, in case, the module) of the endomorphisms $\Hom_\QQ(\MM,\MM)$ will be denoted by $\Edm(\MM) = \la \edm(\MM), \sqcup, \bot^\bot \ra$.
\end{definition}

By Proposition~\ref{dualhomomq}, the sup-lattices $\Hom_\QQ(\MM,\NN)$ and $\Hom_\QQ(\NN^{\op},\MM^{\op})$ are isomorphic and, if $\QQ$ is a commutative quantale, they are isomorphic $\QQ$-modules.

\begin{proposition}\label{coprodmq}
Let $\QQ$ be a quantale. For any family of $\QQ$-modules $\{\MM_i\}_{i \in I}$, the coproduct $\coprod_{i \in I} \MM_i$ is the product $\prod_{i \in I} \MM_i$ equipped with the right inverses $\mu_i: \MM_i \lto \prod_{i \in I} \MM_i$ of the projections $\pi_i: \prod_{i \in I} \MM_i \lto \MM_i$. Moreover, for all $i \in I$, $\pi_i \circ \mu_i = \id_{\MM_i}$.
\end{proposition}
\begin{proof}
First of all, let us observe that, for any fixed index $j \in I$ and for any $x \in M_j$, $\mu_j(x)$ is the family of $\prod_{i \in I} \MM_i$ in which all the elements are equal to $\bot$ except the $j$-th that is equal to $x$. Thus, in particular, for any family $(x_i)_{i \in I}$, $(x_i)_{i \in I} = \bigvee_{i \in I} \mu_i(x_i)$.

Now we want to prove that, given an arbitrary $\QQ$-module $\MM$ and a family of homomorphisms $f_i: \MM_i \lto \MM$, there exists a unique homomorphism $f: \prod_{i \in I} \MM_i \lto \MM$ such that $f \circ \mu_i = f_i$ for all $i \in I$.

Let $f\left((x_i)_{i \in I}\right) = \bigvee_{i \in I} f_i(x_i)$. For all $i \in I$ and $x_i \in M_i$, we have $(f \circ \mu_i)(x_i) = f\left(\mu_i(x_i)\right) = f_i(x_i) \vee \bigvee_{j \in I \setminus \{i\}} f_j(\bot) = f_i(x_i)$.

Now, let $f': \prod_{i \in I} \MM_i \lto \MM$ be another homomorphism such that $f' \circ \mu_i = f_i$ for all $i \in I$. Then $f'\left((x_i)_{i \in I}\right) = f'\left(\bigvee_{i \in I} \mu_i(x_i)\right) = \bigvee_{i \in I} f'(\mu_i(x_i)) = \bigvee_{i \in I} f_i(x_i) = f\left((x_i)_{i \in I}\right)$. Thus $f' = f$ and also the uniqueness is proved.
\end{proof}

\section{$\Q$-module structural closure operators and transforms}
\label{nucleisec}

In Section~\ref{resmapsec} we saw the connection between the homomorphisms whose domain is a given sup-lattice $\LL$ and the closure operators on $\LL$. Now, since $\Q$-module homomorphisms are special sup-lattice homomorphisms, we expect $\Q$-modules to be fit for a suitable class of operators that would allow us to establish a similar result. Indeed the notion that corresponds to that of closure operator, in the categories of $\Q$-modules, is the one of \emph{structural closure operator}, also called a \emph{nucleus} (see~\cite{paseka}).

As a matter of fact, the importance of such operators goes far beyond the relationship with homomorphisms. In~\cite{galatostsinakis} the authors proved that structural closure operators between complete posets allow an algebraic representation of propositional deductive systems, and a re-examination of these results in the framework of $\Q$-modules has been presented in~\cite{russo}. Moreover, they are also involved in the applications to image processing we will present in Section~\ref{appsec}.

\begin{definition}\label{structural}
Let $\QQ$ be a quantale and $\MM$ a $\QQ$-module. A map $\g: M \lto M$ is said to be a \emph{structural operator} on $\MM$ provided it satisfies, for all $m, n \in M$ and $q \in Q$, the following conditions:
\begin{enumerate}
\item[$(i)$]$m \leq \g(m)$;
\item[$(ii)$]$m \leq n$ implies $\g(m) \leq \g(n)$;
\item[$(iii)$]$q \star \g(m) \leq \g(q \star m)$.
\end{enumerate}
We will say that a structural operator $\g$ is a \emph{nucleus}, or a \emph{structural closure operator}, if it is also idempotent:
\begin{enumerate}
\item[$(iv)$]$\g \circ \g = \g$.
\end{enumerate}
Then a nucleus is a closure operator that satisfies also condition $(iii)$. So it is natural to call a \emph{conucleus} a \emph{structural coclosure operator}, i.e. a coclosure operator satisfying condition $(iii)$.
\end{definition}
If $\g$ is a nucleus, we denote by $M_\g$ the $\g$-closed system $\g[M]$ and it is easily seen that $M_\g$ is closed under arbitrary meets. Dually, the image $M^\d$ of a conucleus $\d$ is closed under arbitrary joins. In the following result we give several characterizations of structural closure operators.
\begin{lemma}\emph{\cite{galatostsinakis}}\label{modclop}
Let $\MM$ be a $\QQ$-module and let $\g$ be a closure operator on $\MM$. The following are equivalent:
\begin{enumerate}
\item[$(a)$]$\g$ is structural;
\item[$(b)$]$\g(q \star \g(m)) = \g (q \star m)$, for all $q \in Q$ and $m \in M$;
\item[$(c)$]$\g(m) \ost n = \g(m) \ost \g(n)$, for all $m, n \in M$;
\item[$(d)$]$\g(q \ust m) \leq q \ust \g(m)$, for all $q \in Q$ and $m \in M$;
\item[$(e)$]$q \ust \g(m) \in M_\g$, for all $q \in Q$ and $m \in M$. 
\end{enumerate} 
\end{lemma}

Before we prove the next result, recall that~--- according to Proposition~\ref{dualhomomq}~--- any $\Q$-module homomorphism is a residuated map.

\begin{theorem}\label{homoandnuclei}
Let $\QQ$ be a quantale, $\MM$ and $\NN$ $\QQ$-modules, and $f: \MM \lto \NN$ a homomorphism. Then $f_* \circ f$ is a nucleus on $\MM$.

Vice versa, if $\g$ is a nucleus on $\MM$, then $M_\g$~--- with the join $\vee_\g = \g \circ \vee$, the external product $\star_\g = \g \circ \star$ and the bottom element $\bot_\g = \g(\bot)$~--- is a $\QQ$-module (denoted by $\MM_\g$) and there exists $f_\g \in \hom_\QQ(\MM, \MM_\g)$ such that $f_{\g\ast} \circ f_\g = \g$.
\end{theorem}
\begin{proof}
By Corollary~\ref{resclosureinterior}, if $(f,f_*)$ is an adjoint pair, then $f_* \circ f$ is a closure operator on the domain of $f$; therefore what we need to prove is condition $(iii)$ of Definition~\ref{structural} for $f_* \circ f$.

Since $(f,f_*)$ is an adjoint pair and $f$ is a homomorphism,
\begin{eqnarray}
&& q \star f(m) \leq q \star f(m) \nonumber \\ 
&\liff \quad& q \star (f \circ f_* \circ f)(m) \leq q \star f(m) \nonumber \\
&\liff \quad& f(q \star (f_* \circ f)(m)) \leq q \star f(m) \nonumber \\
&\liff \quad& q \star (f_* \circ f)(m) \leq f_* (q \star f(m)) \nonumber \\
&\liff \quad& q \star (f_* \circ f)(m) \leq (f_* \circ f)(q \star m); \nonumber \\ \nonumber
\end{eqnarray}
then $f_* \circ f$ is a nucleus.

Now let $\g$ be a nucleus on $\MM$. The fact that $\MM_\g$ is a $\QQ$-module comes easily from the fact that it is closed under arbitrary meets, and the definitions of join and product. For the same reasons the map
\begin{equation*}
f_\g: m \in M \lmapsto \g(m) \in M_\g
\end{equation*}
is a $\QQ$-module homomorphism. If we apply Theorem~\ref{residchar}$(b)$ to $f_\g$, we get
$$f_{\g\ast}(\g(m)) = \bigvee \{n \in M \mid f_\g(n) = \g(n) \leq \g(m)\} = \g(m),$$
for any element $\g(m) \in M_\g$; thus $f_{\g\ast} = \id_{M \restr M_\g}$, whence
$$(f_{\g\ast} \circ f_\g)(m) = (\id_{M \restr M_\g} \circ f_\g)(m) = \id_{M \restr M_\g}(\g(m)) = \g(m),$$
for all $m \in M$. The theorem is proved.
\end{proof}

In the same hypotheses of the previous theorem, we observe explicitly that, even if $f \circ f_*$ is an interior operator in the sup-lattice $\NN$, it is not~--- in general~--- a conucleus on the $\QQ$-module $\NN$. It is, instead, a nucleus on the $\QQ$-module $\NN^{\op}$. It suffices to notice that, by Proposition~\ref{dualhomomq}, $f_*$ is a $\QQ$-module homomorphism from $\NN^{\op}$ to $\MM^{\op}$ whose residual is $f$, and to apply Theorem~\ref{homoandnuclei}.

The next result is an interesting, though immediate, property of structural closure operators.
\begin{proposition}\emph{\cite{galatostsinakis}}\label{gammacyclic}
Let $\QQ$ be a quantale. If $\g: \QQ \lto \QQ$ is a nucleus on the $\QQ$-module $\QQ$, then $\QQ_\g$ is a cyclic module and it is generated by $\g(e)$.
\end{proposition}

Next, we introduce the $\Q$-module transforms and prove some results about them. Then we show that any direct transform is a $\Q$-module homomorphism and a residuated map whose residual is its inverse transform.

On the other hand, $\Q$-module \emph{faithful} transforms~--- i.e. those transforms whose kernel is a coder (see Definition~\ref{coders})~--- have many further interesting properties. We will also see that the \L ukasiewicz transform defined in~\cite{dinolarusso} is effectively a (orthonormal) $\Q$-module transform.

\begin{definition}\label{qwtransform}
Let $\QQ \in \Q$ and $X, Y$ be non-empty sets and let us consider the free $\QQ$-modules $\QQ^X$ and $\QQ^Y$. We call a \emph{$\Q$-module transform} between $\QQ^X$ and $\QQ^Y$, with \emph{kernel} $p$, the map
$$H_p: Q^X \lto Q^Y$$
defined by
\begin{equation}\label{qwtransformeq}
H_p f(y) = \bigvee_{x \in X} f(x) \cdot p(x,y) \quad \textrm{for all } y \in Y,
\end{equation}
where $p \in Q^{X \times Y}$. Its \emph{inverse transform} $\Lambda_p: Q^Y \lto Q^X$ is the map defined by
\begin{equation}\label{qwinverse transformeq}
\Lambda_p g(x) = \bigwedge_{y \in Y} g(y) / p(x,y) \quad \textrm{for all } x \in X.
\end{equation}
\end{definition}

\begin{remark}\label{righttrans}
Recalling that we are using the notations of left modules, we observe that, if we consider $\QQ^X$ and $\QQ^Y$ as right modules, the direct and inverse transforms are defined respectively by
\begin{equation}\label{qwtransformeqright}
H_p f(y) = \bigvee_{x \in X} p(x,y) \cdot f(x) \quad \textrm{for all } y \in Y,
\end{equation}
and
\begin{equation}\label{qwinverse transformeqright}
\Lambda_p g(x) = \bigwedge_{y \in Y} p(x,y) \under g(y) \quad \textrm{for all } x \in X.
\end{equation}
Up to a suitable reformulation, all the results we will present for $\Q$-module transforms hold also for free right modules.
\end{remark}

\begin{theorem}\label{wadjointpair}
Let $\QQ \in \Q$, $X, Y$ be two non-empty sets and $p \in Q^{X \times Y}$. If $H_p$ is the $\Q$-module transform, with kernel $p$, between $\QQ^X$ and $\QQ^Y$, and $\Lambda_p$ is its inverse transform, then the following hold:
\begin{enumerate}
\item[$(i)$]$(H_p,\Lambda_p)$ is an adjoint pair, i.e. $H_p$ is a residuated map and $\Lambda_p = H_{p\ast}$;
\item[$(ii)$]$H_p \in \hom_{\lMQ}\left(\QQ^X,\QQ^Y\right)$ and $\Lambda_p \in \hom_{\rMQ}\left(\left(\QQ^Y\right)^{\op},\left(\QQ^X\right)^{\op}\right)$;
\item[$(iii)$]$\Lambda_p \circ H_p$ is a nucleus over $\QQ^X$ and $H_p \circ \Lambda_p$ is a nucleus over $\left(\QQ^Y\right)^{\op}$.
\end{enumerate} 
\end{theorem}
\begin{proof}
\begin{enumerate}
\item[$(i)$] Both $H_p$ and $\Lambda_p$ are clearly monotone; let us prove that $(H_p,\Lambda_p)$ is an adjoint pair by showing that (\ref{resid1}) and (\ref{resid2}) hold. For any $f \in Q^X$ and $g \in Q^Y$, we have:
$$\begin{array}{lll}
H_p f \leq g && \iff \quad \\
H_p f(y) \leq g(y) & \forall y \in Y & \iff \\
\bigvee_{x \in X} f(x) \cdot p_Y(x,y) \leq g(y) & \forall y \in Y & \iff  \\
f(x) \cdot p_Y(x,y) \leq g(y) & \forall x \in X, \forall y \in Y & \iff  \\
f(x) \leq g(y)/p_Y(x,y) & \forall x \in X, \forall y \in Y & \iff \\
f(x) \leq \bigwedge_{y \in Y} g(y)/p_Y(x,y) = \Lambda_p g(x) & \forall x \in X. &
\end{array}$$
Hence $H_p f \leq g \iff f \leq \Lambda_p g$; thus, by setting alternatively $g = H_p f$ and $f = \Lambda_p g$, we get respectively (\ref{resid2}) and (\ref{resid1}), and the $(i)$ is proved.
\item[$(ii)$] Since $H_p$ is a residuated map, it is a sup-lattice homomorphism. Moreover it is evident from the definition that $H_p$ preserves the scalar multiplication.
\item[$(iii)$] It follows from the $(ii)$ and Theorem~\ref{homoandnuclei}.
\end{enumerate}
\end{proof}

The following classification of the kernels has few interesting theoretical implications but it is important for applications to image processing.
\begin{definition}\label{coders}
Let $\QQ \in \Q$, and $X, Y$ be non-empty sets. Let us consider a map $p \in Q^{X \times Y}$; we set the following definitions:
\begin{enumerate}
\item[$(i)$]$p$ is called a \emph{coder} iff there exists an injective map $\e: Y \lto X$ such that $e \leq p(\e(y),y)$ for all $y \in Y$;
\item[$(ii)$]$p$ is said to be \emph{normal} iff there exists an injective map $\e: Y \lto X$ such that $p(\e(y),y) = e$ for all $y \in Y$;
\item[$(iii)$]$p$ is said to be \emph{strong} iff it is normal and
\begin{equation}\label{strong}
p(\e(y_1),y_2) = \bot \quad \textrm{for all } y_1, y_2 \in Y, y_1 \neq y_2;
\end{equation}
\item[$(iv)$]$p$ is said to be \emph{orthogonal} iff $p(x,y_1) \cdot p(x,y_2) = \bot$ for all $y_1, y_2 \in Y$ such that $y_1 \neq y_2$ and for all $x \in X$;
\item[$(v)$]$p$ is said to be \emph{orthonormal} iff it is orthogonal and normal.
\end{enumerate}
If $p$ is a coder, the $\Q$-module transform $H_p$ will be called \emph{faithful}.
\end{definition}

\begin{remark}\label{orthonormal}
\begin{enumerate}
\item[$(i)$] If $p$ is normal, then it is a coder.
\item[$(ii)$] If $p$ is strong, it is a normal coder by definition.
\item[$(iii)$]If $p$ is an orthonormal map and $\e: Y \lto X$ is an injective map as in Definition~\ref{coders} $(ii)$, for any two arbitrary different elements $y_1, y_2 \in Y$, from $p(\e(y_1),y_1) \cdot p(\e(y_1),y_2) = e \cdot p(\e(y_1),y_2) = \bot$, it follows that $p(\e(y_1),y_2) = \bot$. Then any orthonormal map is a strong coder.
\end{enumerate}
\end{remark}

\begin{definition}\label{qtransform}
Let $\QQ \in \Q$, $X, Y$ be two non-empty sets and $p \in Q^{X \times Y}$ be a coder. Let us consider the faithful $\Q$-module transform $H_p$, with kernel $p$, between $Q^X$ and $Q^Y$. We set the following definitions:
\begin{enumerate}
\item[$(i)$]$H_p$ will be called a \emph{normal transform} if $p$ is normal; 
\item[$(ii)$]$H_p$ will be called a \emph{strong transform} if $p$ is strong; 
\item[$(iii)$]$H_p$ will be called an \emph{orthonormal transform} if $p$ is orthonormal.
\end{enumerate}
\end{definition}

\begin{exm}\label{luktransf}
In~\cite{dinolarusso} the authors define the \emph{\L ukasiewicz transform}
$$H_n: [0,1]^{[0,1]} \lto [0,1]^n$$
and its inverse transform $\Lambda_n$, respectively, as
$$H_n(f) = \left(\bigvee_{x \in [0,1]} f(x) \odot p_0(x), \ldots, \bigvee_{x \in [0,1]} f(x) \odot p_{n-1}(x)\right),$$
$$\Lambda_n(v_0, \ldots, v_n) = \bigwedge_{k=0}^{n-1} p_k(x) \to v_k = \left(\bigvee_{k=0}^{n-1} v_k^* \odot p_k(x)\right)^*,$$
where
\begin{equation*}
p_{0}(x)= \begin{cases} -(n-1)x+1 & \text{if $0\leq x\leq \frac{1}{n-1}$ } \\
     0 & \text{otherwise}
\end{cases},
\end{equation*}

\begin{equation*}
p_{n-1}(x)=\begin{cases}

(n-1)x-(n-2)& \text{if $\frac{n-2}{n-1}\leq x\leq 1$ } \\

0 & \text{otherwise}
\end{cases}
\end{equation*}
and, for $k=1,\ldots, n-2$,
\begin{equation*}
p_{k}(x)=\begin{cases}
     (n-1)x-(k-1) & \text{if  $\frac{k-1}{n-1}\leq x\leq   \frac{k}{n-1}$ }
\\
      -(n-1)x+k+1 & \text{if  $\frac{k}{n-1}\leq x\leq \frac{k+1}{n-1} $} \\
 0 & \text{otherwise}
\end{cases}.
\end{equation*}

It is easy to see that $\la [0,1]^{[0,1]}, \vee, \0\ra$ and $\la [0,1]^n, \vee, \0\ra$ are $\Q$-modules over the commutative quantale $\la [0,1], \vee, \odot, 0, 1 \ra$ and $H_n$ is a $\Q$-module orthonormal transform with inverse transform $\Lambda_n$, if we set:
\begin{enumerate}
\item[-] $\I_n = \{0, \ldots, n-1\}$,
\item[-] $\e: k \in \I_n  \lmapsto \frac{k}{n-1} \in [0,1]$,
\item[-] $p_{I_n}(x,k) = p_k(x)$ for all $x \in [0,1]$ and for all $k \in \{0, \ldots, n-1\}$.
\end{enumerate}
\end{exm}

\begin{theorem}\label{sadjointpair}
Let $\QQ \in \Q$ and let $H_p$ be a $\Q$-module strong transform, by the coder $p \in Q^{X \times Y}$, with inverse transform $\Lambda_p$. Then
\begin{equation*}
H_p \circ \Lambda_p = \id_{Q^Y};
\end{equation*}
thus $H_p$ is onto and, by Proposition~\ref{residprop} $(iii)$, $\Lambda_p$ is one-one.
\end{theorem}
\begin{proof}
By Theorem~\ref{wadjointpair} we have $H_p \circ \Lambda_p \leq \id_{Q^Y}$. In order to prove the inverse inequality, let us proceed as follows.

Since $p$ is strong, we can consider an injective map $\e: Y \lto X$ such that condition (\ref{strong}) holds. Let now $g \in Q^Y$ be an arbitrary function and let us fix an arbitrary $\ov y \in Y$. We have:
\begin{eqnarray}
\lefteqn{(H_p \circ \Lambda_p) g(\ov y) =} \nonumber \\
&& \bigvee_{x \in X}\left(\left(\bigwedge_{y \in Y} g(y) / p(x,y)\right) \cdot p(x,\ov y)\right) \geq \nonumber \\
&& \left(\bigwedge_{y \in Y} g(y) / p(\e(\ov y),y)\right) \cdot p(\e(\ov y),\ov y) = \nonumber \\
&& \left(\bigwedge_{\mathop{}^{y \in Y}_{y \neq \ov y}} g(y) / p(\e(\ov y),y)\right) \wedge \left(g(\ov y) / p(\e(\ov y),\ov y)\right) = \nonumber \\
&& \left(\bigwedge_{\mathop{}^{y \in Y}_{y \neq \ov y}} g(y) / \bot\right) \wedge \left(g(\ov y) / e\right) = \nonumber \\
&& \top \wedge g(\ov y) = \nonumber \\
&& g(\ov y) \nonumber.
\end{eqnarray}
Since the above relations hold for all $g \in Q^Y$ and $\ov{y} \in Y$, the result is proved.
\end{proof}

\begin{lemma}\label{transf=coder}
Let $\QQ \in \Q$, $X$ and $Y$ be non-empty sets and $p, p' \in Q^{X \times Y}$ be two maps. Then $H_p = H_{p'}$ if and only if $p = p'$.
\end{lemma}
\begin{proof}
Since the other implication is trivial, let us prove that $H_p = H_{p'}$ implies $p = p'$ by showing that, if $p \neq p'$, then $H_p \neq H_{p'}$.

By assumption, there exists a pair $(\ov x,\ov y) \in X \times Y$ such that $p(\ov x, \ov y) \neq p'(\ov x, \ov y)$. Let us consider the map $f \in Q^X$ defined by
$$f(x) = \left\{
			\begin{array}{ll}
			e & \textrm{if } x = \ov x \\
			\bot & \textrm{if } x \in X \setminus \{\ov x\} \\
			\end{array}.
			\right.
$$
It is immediate to verify that $H_p f(\ov y) = p(\ov x, \ov y) \neq p'(\ov x, \ov y) = H_{p'} f(\ov y)$, and the result follows.
\end{proof}

The previous result ensures that a $\Q$-module transform $H_p$ is completely determined by its kernel $p$.

In what follows we will always assume that $Y$ is a subset of $X$ and, if $p \in Q^{X \times Y}$ is a coder, then the map $\e$ is the inclusion map $\id_{X \restr Y}: y \in Y \lmapsto y \in X$.

\begin{lemma}\label{projection}
Let $\QQ \in \Q$, $X$ be a non-empty set, $Y$ be a non-empty subset of $X$ and $p \in Q^{X \times Y}$ be a coder. Then, for any fixed $y \in Y$, $H_p f(y) = f(y)$ for all $f \in Q^X$ if and only if
\begin{equation}\label{projcoder}
p(x,y) = \left\{
				\begin{array}{ll}
				e & \textrm{if } x = y \\
				\bot & \textrm{if } x \neq y
				\end{array}.\right.
\end{equation}
\end{lemma}
\begin{proof}
If condition (\ref{projcoder}) holds, then $H_p f(y) = f(y)$ for all $f \in Q^X$, trivially.

On the other hand, if (\ref{projcoder}) does not hold for $p$, then we distinguish two cases:
\begin{enumerate}
\item[Case 1:]$p(y, y) = q_1 \neq e$, for some $q_1 \in Q$;
\item[Case 2:]there exists $\ov x \in X$ such that $p(\ov x, y) = q_2 \gneqq \bot$.
\end{enumerate}
In the first case, let $f \in Q^X$ be the map defined by
$$f(x) = \left\{
			\begin{array}{ll}
			e & \textrm{if } x = y \\
			\bot & \textrm{if } x \neq y \\
			\end{array}.
			\right.
$$
It is easy to see that $H_p f(y) = q_1 \neq e = f(y)$.

In the second case, let $g \in Q^X$ be the map defined as
$$g(x) = \left\{
			\begin{array}{ll}
			e & \textrm{if } x = \ov x \\
			\bot & \textrm{if } x \neq \ov x \\
			\end{array}.
			\right.
$$
Clearly $H_p g(y) = q_2 \neq \bot = g(y)$, and the lemma is proved.
\end{proof}

Given the sets $X$ and $Y$, we will denote by $\pi_Y$ the coder defined, for all $x \in X$ and for all $y \in Y$, by condition (\ref{projcoder}) and we will call it a \emph{projective coder}. By Lemma~\ref{projection}, for all $f \in Q^X$, $H_{\pi_Y} f = f_{\restr Y}$, i.e. $H_{\pi_Y}$ is the projection of $Q^X$ on $Q^Y$.

\begin{definition}\label{psupport}
If $p$ is a coder of $Q^{X \times Y}$, let $Y_p' \subseteq Y$ be the set of all the elements $y$ of $Y$ such that $p(x,y)$ is defined by condition (\ref{projcoder}):
$$Y_p' = \left\{y \in Y \ | \ p(x,y) = \pi_Y(x,y) \right\}.$$
The set $\dot Y_p = Y \setminus Y_p'$ will be called the \emph{support} of $p$ and the restriction $\ul p = p_{\restr \dot Y_p}$ will be called the \emph{core} of $p$. If $\dot Y_p = Y$, then $p = \ul p$ and we will say that $p$ is \emph{irreducible}; $p$ is \emph{reducible} if $\dot Y_p \subsetneqq Y$.
\end{definition}

\begin{definition}\label{closurecod}
Given a coder $p \in Q^{X \times Y}$ and a set $Z$ such that $Y \subseteq Z \subseteq X$, let us consider the extension $p^Z$ of the coder $p$ to $X \times Z$, defined as follows:
\begin{equation*}
p^Z(x,z) = \left\{
						\begin{array}{ll}
							p(x,z) & \textrm{if } (x,z) \in X \times Y \\
							\pi_{Z \setminus Y}(x,z) & \textrm{if } (x,z) \in (X \times Z) \setminus (X \times Y) \\
						\end{array}.
						\right.
\end{equation*}
The coder $p^Z$ will be called the \emph{projective extension of $p$ to $Z$}. In this case, it is clear that $\dot Y_p = \dot Z_{p^Z}$ and $\ul p = \ul p^Z$.

If $Z = X$, we will denote $p^X$ by $\ov p$ and we will call it the \emph{closure} of $p$. So $\ov p$ is the coder defined by 
\begin{equation*}
\ov p(x,z) = \left\{
						\begin{array}{ll}
							p(x,z) & \textrm{if } (x,z) \in X \times Y \\
							\pi_{X \setminus Y}(x,z) & \textrm{if } (x,z) \in (X \times X) \setminus (X \times Y) \\
						\end{array}.
						\right.
\end{equation*}
Clearly, for any coder $p \in Q^{X \times X}$, $\ov p = p$; therefore, such coders will be called \emph{closed coders}.
\end{definition}

\begin{definition}\label{equiproj}
Let $\QQ \in \Q$, $X$ be a non-empty set, $Y, Z$ be two non-empty subsets of $X$ and $p \in Q^{X \times Y}$, $p' \in Q^{X \times Z}$ be two coders. We will say that $p$ and $p'$ are \emph{equivalent up to projections}~--- and we will write $p \doteq p'$~--- iff $\ul p = \ul p'$, i.e. iff $\dot Y_p = \dot Z_{p'}$ and $p_{\restr \dot Y_p} = p'_{\restr \dot Y_p}$.
\end{definition}

\begin{proposition}\label{eqclosure}
Let $\QQ \in \Q$, $X$ be a non-empty set, $Y, Z$ be two non-empty subsets of $X$ and $p \in Q^{X \times Y}$, $p' \in Q^Z$ be two coders. Then
$$p \doteq p' \quad \iff \quad \ov p = \ov p'.$$
In other words, $p$ and $p'$ are equivalent up to projections if and only if they have the same closure.
\end{proposition}
\begin{proof}
It is trivial.
\end{proof}
The last definitions and Proposition~\ref{eqclosure} are significant, again, for applications. In the next result we invert Theorem~\ref{wadjointpair}, showing that all the homomorphisms between free modules are transforms.

\begin{theorem}\label{repr}
The sup-lattices $\Hom_\QQ(\QQ^X,\QQ^Y)$ and $\QQ^{X \times Y}$ are isomorphic.
\end{theorem}
\begin{proof}
Let
\begin{equation}\label{hbarxy}
\hbar: \ Q^{X \times Y} \ \lto \ \hom_\QQ(\QQ^X,\QQ^Y)
\end{equation}
be the map defined by $\hbar(p) = H_p$, for all $p \in Q^{X \times Y}$; in other words $\hbar$ sends every map $p \in Q^{X \times Y}$ in the transform between $\QQ^X$ and $\QQ^Y$ whose kernel~is~$p$.

The fact that $\hbar$ is injective comes directly from Lemma~\ref{transf=coder}. Moreover it is clear that $\hbar(\bot^{X \times Y}) = \bot^\bot$. Now let $\{k_i\}_{i \in I} \subseteq Q^{X \times Y}$; we must prove that $\hbar\left(\bigvee_{i \in I} k_i \right) = \bigsqcup_{i \in I}\hbar(k_i)$. For any $f \in Q^X$ and for all $y \in Y$, we have
\begin{eqnarray}
\lefteqn{\hbar\left(\bigvee_{i \in I} k_i \right) f(y) = \bigvee_{x \in X} f(x) \cdot \left(\bigvee_{i \in I} k_i \right)(x,y)} \nonumber\\
&&= \bigvee_{x \in X} f(x) \cdot \left(\bigvee_{i \in I} k_i(x,y) \right) = \bigvee_{x \in X} \bigvee_{i \in I} f(x) \cdot k_i(x,y) \nonumber\\
&&= \bigvee_{i \in I} \bigvee_{x \in X} f(x) \cdot k_i(x,y) = \bigvee_{i \in I} \hbar(k_i)f(y) \nonumber\\
&&= \left(\bigvee_{i \in I} \hbar(k_i)f\right)(y) = \left(\bigsqcup_{i \in I} \hbar(k_i)\right)f(y); \nonumber
\end{eqnarray}
whence $\hbar$ is a sup-lattice monomorphism. 

Now we must prove that $\hbar$ is surjective too. Let $h \in \hom_\QQ(\QQ^X,\QQ^Y)$ and, for any $x \in X$, let us consider the map $\chi_x$ defined by (\ref{chi}). Let now $k^h \in Q^{X \times Y}$ be the function defined by
\begin{equation*}
k^h(x,y) = h \chi_x(y), \qquad \textrm{for all } (x,y) \in Q^{X \times Y};
\end{equation*}
then we have
\begin{eqnarray}
\lefteqn{h f(y) = h\left(\bigvee_{x \in X} f(x) \star \chi_x\right)(y)} \nonumber \\
&& = \left(\bigvee_{x \in X} f(x) \star h \chi_x\right)(y) = \bigvee_{x \in X} f(x) \cdot h \chi_x(y) \nonumber \\
&& = \bigvee_{x \in X} f(x) \cdot k^h(x,y) = H_{k^h} f(y), \nonumber \\ \nonumber
\end{eqnarray}
for all $f \in Q^X$ and for all $y \in Y$. It follows that $h = H_{k^h} = \hbar(k^h)$, hence $\hbar$ is a sup-lattice isomorphism whose inverse is~--- obviously~--- defined by 
$$\hbar^{-1} h(x,y) = h \chi_x(y),$$
for all $h \in \hom_\QQ(\QQ^X,\QQ^Y)$ and $(x,y) \in X \times Y$.
\end{proof}

The previous theorem allows us to define the structure of a $\QQ$-module on the sup-lattice $\hom_\QQ(\QQ^X,\QQ^Y)$, also when $\QQ$ is not commutative, by defining the external multiplication
\begin{equation}\label{asthom}
\begin{array}{cccc}
\star: & Q \times \hom_\QQ(\QQ^X,\QQ^Y) & \lto 		& \hom_\QQ(\QQ^X,\QQ^Y) \\
			& (q,h) 												 & \lmapsto & \hbar(q \star \hbar^{-1}(h))
\end{array},
\end{equation}
in such a way that this $\QQ$-module is isomorphic to $\QQ^{X \times Y}$; we will denote this structure by $\Hom_\QQ^\star(\QQ^X,\QQ^Y)$.

We will not discuss projective objects in the categories of quantale modules, here, several interesting results about them were presented in \cite{galatostsinakis}. However, before extending Theorem \ref{repr}, we recall that free modules are projective and that every module is the homomorphic image of a free module.

\begin{theorem}\label{projrepr}
Let $\MM$ and $\NN$ be $\QQ$-modules, $X$ and $Y$ be two sets such that $\MM$ and $\NN$ are homomorphic images of $\QQ^X$ and $\QQ^Y$ respectively, and $\pi: \QQ^X \lto \MM$ and $\pi': \QQ^Y \lto \NN$ be the respective surjective morphisms.

Then, for any homomorphism $h: \MM \lto \NN$ there exists $k \in \QQ^{X \times Y}$ such that $h \circ \pi = \pi' \circ H_k$, where $H_k$ is the transform from $\QQ^X$ to $\QQ^Y$ whose kernel is $k$.
\end{theorem}
\begin{proof}
Consider the following diagram
\begin{equation*}
\xymatrix{
\QQ^X \ar@{-->}[rr]^{H_k} \ar[dd]_\pi && \QQ^Y \ar[dd]^{\pi'} \\
&&\\
\MM	\ar[rr]_h && \NN
}.
\end{equation*}
The existence of the morphism $H_k$ that closes such a diagram follows immediately from the projectivity of $\QQ^X$, and we know by Theorem~\ref{repr} that $H_k$ is indeed a transform.

Nonetheless it is interesting to notice that, since each element $m$ of $\MM$ can be written as $\bigvee_{x \in X} q_x \star \pi(\chi_x)$, with the $q_x$'s in $\QQ$, then
$$\begin{array}{lllll}
h(m) & = & h\left(\bigvee_{x \in X} q_x \star \pi(\chi_x)\right) & = & \bigvee_{x \in X} q_x \star h(\pi(\chi_x)) \\
		 & = & \bigvee_{x \in X} q_x \star \pi'(H_k(\chi_x)) & = & \bigvee_{x \in X} q_x \star \pi'(k(x,{}_-)).
\end{array}$$
This shows how ``concretely'' $k$ determines $h$.
\end{proof}

\section{Applications}
\label{appsec}

In this section we show how certain techniques of image processing, with different scopes, can be grouped together under the common ``algebraic roof'' of $\Q$-module transforms.

The theory of \emph{fuzzy relation equations}~\cite{dinolasessa}, is involved in many algorithms for compression and reconstruction of digital images (see, for example,~\cite{hirotanobuhara,nobuharapedrycz1,hirotapedrycz}). As a matter of fact, fuzzy relations fit the problem of processing the representation of an image as a matrix with the range of its elements previously normalized to $[0,1]$. In such techniques, however, the approach is mainly experimental and the algebraic context is seldom clearly defined.

A first unification of fuzzy image processing has been proposed by I. Perfilieva in~\cite{perfilieva}, with an approach that is analytical rather than algebraic. Moreover, the field of applications of the operators (called \emph{Fuzzy transforms}) defined in~\cite{perfilieva} is limited to the real unit interval, $[0,1]$, endowed with the usual order relation and a triangular norm.

Basically, most of the fuzzy algorithms for image compression, make use of join-product operators, and they can be seen as approximate discrete solutions of fuzzy relation equations of the form $A(x,z) = \bigvee_y B(x,y) \cdot C(y,z)$; after all, a complete lattice order and a multiplication that is residuated w.r.t. the lattice-order are the fundamental ingredients of these operators. So it is natural to think of them as examples of $\Q$-module transforms. Indeed we will see in Subection~\ref{imagesec} that the class of $\Q$-module transforms contains all these operators and much more. 

Further classes of operators that fall within $\Q$-module transforms are those of \emph{mathematical morphological operators}. Mathematical morphology is a technique for image processing and analysis whose origins can be traced back to the book~\cite{matheron}, of 1975, by G. Matheron, and whose development is due mainly to the works by J. Serra and H. J. A. M. Heijmans.

Essentially, mathematical morphological operators analyse the objects in an image by ``probing'' them with a small geometric ``model-shape'' (e.g., line segment, disc, square) called the \emph{structuring element}. These operators are defined on spaces having both a complete lattice order (set inclusion, in concrete applications) and an external action from another ordered structure (the set of translations); they are also usually coupled in adjoint pairs. A description of such operators in terms of $\Q$-module transforms can easily be anticipated.

In the next subsections, rather than dwelling upon technical details, we will try to give the basic ideas of how fuzzy transforms and mathematical morphological operators work (Subsections~\ref{imagesec} and~\ref{matmorsec}), and then to show~---~in Subsection~\ref{ordapprsec}~---~how $\Q$-module transforms suffice to describe all those techniques.

\subsection{Image Compression and Reconstruction}
\label{imagesec}

In the literature of image compression, the fuzzy approach is based essentially on the theory of fuzzy relation equations, deeply investigated by A. Di Nola, S. Sessa, W. Pedrycz and E. Sanchez in~\cite{dinolasessa}. The underlying idea is very easy: a grey-scale image is basically a matrix in which every element represents a pixel and its value, included in the set $\{0, \ldots, 255\}$ in the case of a 256-bit encoding, is the ``grey-level'', where $0$ corresponds to black, $255$ to white and the other levels are, obviously, as lighter as they are closer to $255$. Then, if we normalize the set $\{0, \ldots, 255\}$ by dividing each element by $255$, grey-scale images can be modeled equivalently as fuzzy relations, fuzzy functions (i.e. $[0,1]$-valued maps) or fuzzy subsets of a given set.

As already mentioned, we will neither cover the wide literature on this subject, nor show how such techniques have been developed in recent years (also because it would be a thankless task). Here we rather want to point out the connection with our work, and the best way to show it is to present the first attempt of unifying all (or most of) these techniques in a common algebraic framework, namely the \emph{fuzzy transforms expressed by residuated lattice operations}, introduced by I. Perfilieva in~\cite{perfilieva}.

A binary operation $\ast: [0,1]^2 \lto [0,1]$ is called a \emph{triangular norm}\index{Triangular norm}\index{T-norm}, \emph{t-norm} for short, provided it verifies the following conditions
\begin{enumerate}
\item[-] commutativity: $x \ast y = y \ast x$;
\item[-] monotonicity: $x \ast y \leq z \ast y$ if $x \leq z$ and $x \ast y \leq x \ast z$ if $y \leq z$;
\item[-] associativity: $x \ast (y \ast z) = (x \ast y) \ast z$;
\item[-] $1$ is the neutral element: $1 \ast x = x = x \ast 1$.
\end{enumerate}
A t-norm $\ast$ is called \emph{left-continuous}\index{Triangular norm!Left-continuous --}\index{T-norm!Left-continuous --} if, for all $\{x_n\}_{n \in \N}, \{y_n\}_{n \in \N} \in [0,1]^\N$,
$$\left(\bigvee_{n \in \N} x_n\right) \ast \left(\bigvee_{n \in \N} y_n\right) = \bigvee_{n \in \N} (x_n \ast y_n).$$
In this case, clearly, $\ast$ is a residuated operation and its residuum (unique, since $\ast$ is commutative) is given by
$$x \to y = \bigvee\{z \in [0,1] \mid z \ast x \leq y\}.$$

Although t-norms are the fuzzy logical analogues of the conjunction of classical logic, here we are mainly interested in them as algebraic operations. The defining conditions of t-norms are exactly the same as those that define a partially ordered (integral) Abelian monoid on the real unit interval $[0,1]$. Therefore some authors call t-norm also the monoidal operation of any partially ordered Abelian monoid; then, in this case, the concept of left-continuity can be substituted by the requirement that the Abelian po-monoid is actually a commutative residuated lattice.

By a \emph{fuzzy partition}\index{Fuzzy!-- partition} of the real unit interval $[0,1]$, we mean an $n$-tuple of fuzzy subsets $A_1, \ldots, A_n$, with $n \geq 2$, identified with their membership functions $A_i: [0,1] \lto [0,1]$ satisfying the following \emph{covering property}
\begin{equation}\label{covering}
\textrm{for all $x \in [0,1]$ there exists $i \leq n$ such that} \quad A_i(x) > 0.
\end{equation}
The membership functions $A_1, \ldots, A_n$ are called the \emph{basic functions} of the partition. There is assumed to exist a finite subset $P \subset [0,1]$, consisting of \emph{nodes} $p_1, \ldots, p_l$ where $l$ is a sufficiently large natural number. Moreover, we assume that $P$ is \emph{sufficiently dense} with respect to the fixed partition, i.e.
\begin{equation}\label{suffdense}
\textrm{for all $i \leq n$ there exists $j \leq l$ such that $A_i(p_j) > 0$}.
\end{equation}
\begin{definition}\label{ftransf}
Let $f \in [0,1]^P$ and $A_1, \ldots, A_n$, $n < l$, be basic functions of a fuzzy partition of $[0,1]$. We say that the $n$-tuple $(F_1^\ua, \ldots, F_n^\ua)$ is the $F^\ua$-transform of $f$ with respect to $A_1, \ldots, A_n$ if, for all $k \leq n$,
\begin{equation}\label{fuptransf}
F^\ua_k = \bigvee_{j=1}^l (A_k(p_j) \ast f(p_j)).
\end{equation}
We say that the $n$-tuple $(F_1^\da, \ldots, F_n^\da)$ is the $F^\da$-transform of $f$ with respect to $A_1, \ldots, A_n$ if, for all $k \leq n$,
\begin{equation}\label{fdowntransf}
F^\da_k = \bigvee_{j=1}^l (A_k(p_j) \to_* f(p_j)).
\end{equation}
\end{definition}
\begin{definition}\label{finvtransf}
Let $f \in [0,1]^P$, $A_1, \ldots, A_n$, with $n < l$, be basic functions of a fuzzy partition of $[0,1]$, and $(F_1^\ua, \ldots, F_n^\ua)$ be the $F^\ua$-transform of $f$ with respect to $A_1, \ldots, A_n$ if, for all $k \leq n$. The map defined, for all $j \leq l$, by
\begin{equation}\label{fupinvtransf}
f^\ua(p_j) = \bigwedge_{k=1}^n (A_k(p_j) \to_* F^\ua_k)
\end{equation}
is called the \emph{inverse $F^\ua$-transform} of $f$.

Let $(F_1^\da, \ldots, F_n^\da)$ be the $F^\da$-transform of $f$ with respect to $A_1, \ldots, A_n$ if, for all $k \leq n$. The map defined, for all $j \leq l$, by
\begin{equation}\label{fdainvtransf}
f^\da(p_j) = \bigvee_{k=1}^n (A_k(p_j) \ast F^\da_k)
\end{equation}
is called the \emph{inverse $F^\da$-transform} of $f$.
\end{definition}

Apart from the definitions above, several results on such tranforms are presented in the cited paper; further algebraic results on join-product composition operators were also presented in the aforementioned paper~\cite{dinolarusso}. We do not list them here since they are essentially special cases of more general results that we presented in Sections~\ref{resmapsec}~and~\ref{nucleisec}.

\subsection{Mathematical Morphology}
\label{matmorsec}

In~\cite{goutsiasheijmans}, the authors state:
\begin{quote}
\small{The basic problem in mathematical morphology is to design nonlinear operators that extract relevant topological or geometric information from images. This requires development of a mathematical model for images and a rigorous theory that describes fundamental properties of the desirable image operators.}
\end{quote}

Then, not surprisingly, images are modeled, in the wake of tradition and intuition, as subspaces or subsets of a suitable space $E$, which is assumed to possess some additional structure (topological space, metric space, graph, etc.), usually depending on the kind of task at hand. We have seen that, in the case of digital image compression, the image space is often modeled as the set of all the functions from a set~--- the set of all the pixels~--- to the real unit interval $[0,1]$. Then, depending on several ``experimental'' factors, the properties of $[0,1]$ involved may be the usual operations, the order relation, t-norms and so on.

In mathematical morphology, the family of binary images is given by $\wp(E)$, where $E$ is, in general, $\R^n$ or $\Z^n$, for some $n \in \N$. In the first case we have continuous binary images, otherwise we are dealing with discrete binary images. The basic relations and operations between images of this type are essentially those between sets, namely set inclusion, union, intersection and so on. As a first example, we can consider an image $X$ that is \emph{hidden} by another image $Y$. Then we can formalize this fact by means of set inclusion: $X \subseteq Y$. Analogously, if we simultaneously consider two images $X$ and $Y$, what we see is their union $X \cup Y$; the \emph{background} of an image $X$ is its complement $X^c$ in the whole space, and the part of an image $Y$ that is not covered by another image $X$ is the set difference $Y \setminus X = Y \cap X^c$.

It is easily anticipated, then, that the lattices are the algebraic structures required for abstracting the ideas introduced so far. Nonetheless, keeping in mind the models $\R^n$ and $\Z^n$, it is possible to introduce the concepts of \emph{translation}\index{Translation!-- of an image} of an image and \emph{translation invariance} of an operator, by means of the algebraic operation of sum.

The reader may recognize the following definitions as those of a residuated map and its residual, and of an adjoint pair.
\begin{definition}\label{dilero}
Let $\LL$, $\MM$ be complete lattices. A map $\d: L \lto M$ is called a \emph{dilation}\index{Dilation} if it distributes over arbitrary joins, i.e., if $\d\left({}^\LL\bigvee_{i \in I} x_i\right) = {}^\MM\bigvee_{i \in I} \d(x_i)$, for every family $\{x_i\}_{i \in I} \subseteq L$. A map $\e: M \lto L$ is called an \emph{erosion}\index{Erosion} if it distributes over arbitrary meets, i.e., if $\e\left({}^\MM\bigwedge_{i \in I} y_i\right) = {}^\LL\bigwedge_{i \in I} \e(y_i)$, for every family $\{y_i\}_{i \in I}$ of elements of $M$.

Two maps $\d: L \lto M$ and $\e: M \lto L$ are said to form an \emph{adjunction}\index{Adjunction}, $(\d,\e)$, between $\LL$ and $\MM$ if $\d(x) \leq y \liff x \leq \e(y)$, for all $x \in L$ and $y \in M$.
\end{definition}
Notice that the notation used in mathematical morphology is slightly different. Indeed, an adjoint pair is presented with the residuated map in the second coordinate and its residual in the first. Here, in order to avoid confusion, we keep on using the notations introduced in Section~\ref{resmapsec}. So we may reformulate the definition above by considering the sup-lattice reducts of $\LL$ and $\MM$, and saying that $\d: L \lto M$ is a dilation if it is a sup-lattice homomorphism between $\LL$ and $\MM$. Dually, an erosion $\e: M \lto L$ is a sup-lattice homomorphism between $\MM^{\op}$ and $\LL^{\op}$. Then a dilation $\d$ and an erosion $\e$ form an adjunction if $\e = \d_*$.

Assume that $\d: \LL \lto \MM$ is a dilation. For $x \in L$, we can write
\begin{equation}\label{dilation}
\d(x) = \bigvee_{y \leq x} \d(y),
\end{equation}
where we have used the fact that $\d$ distributes over join. Every dilation defined on $\LL$ is of the form (\ref{dilation}), and the adjoint erosion is
given by
\begin{equation}\label{erosion}
\e(y) = \bigvee_{\d(x) \leq y} x.
\end{equation}
In the case of powersets, if $\d$ is a dilation between $\wp(E)$ and $\wp(F)$, where $E$ and $F$ are nonempty sets. For $X \subseteq E$, we can write
\begin{equation}\label{wpdilation}
\d(X) = \bigcup_{x \in X} \d(\{x\}),
\end{equation}
and the adjoint erosion is, for all $Y \subseteq F$,
\begin{equation}\label{wperosion}
\e(Y) = \{x \in E \mid \d(\{x\}) \subseteq Y\} = \bigcup_{\d(X) \subseteq Y} X.
\end{equation}

Next, we consider the special case when the operators are translation invariant. In this case, the sets $\d(\{x\})$ are translations of a fixed set, called the \emph{structuring element}\index{Structuring element}, by $\{x\}$. Let $E$ be $\R^n$ or $\Z^n$, and consider the complete lattice $\wp(E)$; given an element $h \in E$, we define the $h$-\emph{translation} $\tau_h$ on $\wp(E)$ by setting, for all $X \in \wp(E)$,
\begin{equation}\label{translation}
\tau_h(X) = X + h = \{x + h \mid x \in X\},
\end{equation}
where the sum is intended to be defined coordinatewise.

An operator $f: \wp(E) \lto \wp(E)$ is called \emph{translation invariant}\index{Translation!-- invariant operator}\index{Operator!Translation invariant --}, \emph{T-invariant} for short, if $\tau_h \circ f = f \circ \tau_h$ for all $h \in E$. It can be proved that every T-invariant dilation on $\wp(E)$ is given by
\begin{equation}\label{tdilation}
\d_A(X) = \bigcup_{x \in X} A + x,
\end{equation}
and every T-invariant erosion is given by
\begin{equation}\label{terosion}
\e_A(X) = \{y \in E \mid A + y \subseteq X\} = \{y \in E \mid y \in X + \breve A\},
\end{equation}
where $A$ is an element of $\wp(E)$, called the \emph{structuring element}, and $\breve A = \{- a \mid a \in A\}$ is the reflection of $A$ around the origin.

Now we observe that the expressions for erosion and dilation in (\ref{tdilation}) and (\ref{terosion}) can also be written, respectively, as
\begin{equation}\label{tdilation2}
\d_A(X)(y) = \bigvee_{x \in E} A(y - x) \wedge X(x)
\end{equation}
and
\begin{equation}\label{terosion2}
\e_A(Y)(x) = \bigwedge_{y \in E} A(y - x) \to Y(y),
\end{equation}
where each subset $X$ of $E$ is identified with its membership function
$$X: x \in E \lmapsto \left\{
\begin{array}{ll}
1 & \textrm{if } x \in X \\
0 & \textrm{if } x \in X^c,
\end{array}\right. \in \{0, 1\}$$
and $X \to Y =: X^c \vee Y$.
Moving from these expressions, and recalling that $\wedge$ is a residuated commutative operation (that is, a continuous t-norm) whose residuum is $\to$, it is possible to extend these operations from the complete lattice of sets $\wp(E) = \{0,1\}^E$ to the complete lattice of fuzzy sets $[0,1]^E$, by means of continuous t-norms and their residua. What we do, concretely, is extend the morphological image operators of dilation and erosion, from the case of binary images, to the case of grey-scale images.

So let $\ast$ be a continuous t-norm and $\to$ be its residuum; a grey-scale image $X$ is a fuzzy subset of $E$, namely a map $X: E \lto [0,1]$. Given a fuzzy subset $A \in [0,1]^E$, called a \emph{fuzzy structuring element}\index{Structuring element!Fuzzy --}\index{Fuzzy!-- structuring element}, the operator
\begin{equation}\label{fdilation}
\d_A(X)(y) = \bigvee_{x \in E} A(y - x) \ast X(x)
\end{equation}
is a translation invariant dilation on $[0,1]^E$, and the operator
\begin{equation}\label{ferosion}
\e_A(X)(x) = \bigwedge_{y \in E} A(y - x) \to X(y)
\end{equation}
is a translation invariant erosion on $[0,1]^E$.

Combining the operators of dilation and erosion by means of the usual algebraic operations in $[0,1]$ it is possible to define new operators, e.g. \emph{outlining} and \emph{top-hat transform}. Their treatment is beyond the scope of this paper, hence we will not present them in details; however some examples can be found at the webpage~\cite{sito3} or in some major reference works in the area of mathematical morphology, such as~\cite{goutsiasheijmans, heijmans1, heijmans, heijmansronse, serra2}, as well as the aforementioned~\cite{serra1}.

\subsection{A unified approach by $\Q$-module transforms}
\label{ordapprsec}

The operators defined so far in this section have a familiar form. Indeed they are all special cases of $\Q$-module transforms between free modules over the quantale reducts of residuated lattice structures defined on the real unit interval $[0,1]$. We now analyse them in detail.

Let us consider the $F^\ua$-transforms of Definition~\ref{ftransf}. Its domain is $[0,1]^l$ and its codomain is $[0,1]^n$ with $n \leq l$. We get immediately that a $\Q$-module transform
$$H_k: f \in [0,1]^l \lmapsto \bigvee_{j = 1}^l f(j) \ast k(j,{}_-) \in [0,1]^n$$
is an $F^\ua$-transform iff the kernel $k$ satisfies condition (\ref{suffdense}) rewritten as
\begin{equation}\label{qsuffdense}
\textrm{for all $i \leq n$ there exists $j \leq l$ such that $k(j,i) > 0$}.
\end{equation}
Obviously, the inverse $F^\ua$-transform of $H_k$ is right
$$\Lambda_k: g \in [0,1]^n \lmapsto \bigwedge_{i = 1}^n k({}_-,i) \to_* g(i) \in [0,1]^l,$$
i.e. the inverse $\Q$-module transform of $H_k$. The case of $F^\da$-transforms is dual to that of $F^\ua$, in the sense that the direct $F^\da$-transform is an inverse $\Q$-module transform, thus a homomorphism between the duals of free modules, and the inverse transform has the shape of a $\Q$-module transform. In other words, for $F^\da$-transforms we assume $l \leq n$ and the condition 
\begin{equation}\label{qsuffdense2}
\textrm{for all $j \leq l$ there exists $i \leq n$ such that $k(j,i) > 0$};
\end{equation}
then $\Lambda_k$ above is the direct $F^\da$-transform and $H_k$ is its inverse.

We already observed in Subsection~\ref{matmorsec} that dilations are precisely the sup-lattice homomorphisms, while erosions are their residua. In order to faithfully represent dilations and erosions that are translation invariant as $\Q$-module transforms from a free $[0,1]$-module to itself, we make the further assumption that the set over which the free module is defined has the additional structure of an Abelian group.

So, let $\mathbf{X} = \la X, +, -, 0 \ra$ be an Abelian group, $\ast$ a t-norm on $[0,1]$, and consider the free $[0,1]$-module $[0,1]^X$. For any element $k \in [0,1]^X$, we define the two variable map $\ov k: (x,y) \in X \times X \lmapsto k(y - x) \in [0,1]$. Then, for all $k \in [0,1]^X$, the translation invariant dilation, on $[0,1]^X$, whose structuring element is $k$, is precisely the $\Q$-module transform $H_{\ov k}$, with the kernel $\ov k$ defined above. Obviously, the translation invariant erosion whose structuring element is $k$ is $\Lambda_{\ov k}$.

Then the representation of both fuzzy transforms and pairs dilation--erosion as quantale module transforms is trivial. Actually, what we want to point out here is that, if we drop the assumption that our quantale is defined on $[0,1]$, the classes of transforms defined in this section become much wider. The purpose of this consideration is not to suggest a purely speculative abstraction but, rather, to underline that suitable generalizations of these operators exist already and they may be useful provided their underlying ideas are extended to other kind of tasks. Indeed the aim of fuzzy transforms is to approximate maps that take values in $[0,1]$; hence the area of application of the whole class of $\Q$-module transform, as approximating operators, can be easily enlarged. On the other hand, the idea of dilating and eroding a shape, in order to analyse it, has not yet found an appropriate concrete extension to situations where $[0,1]$ must be replaced by a non-integral quantale. Nonetheless, we strongly believe (and we are working in this direction) that $\Q$-module dilations and erosions will soon find concrete tasks for being fruitfully applied.

\section*{Conclusion}

In this paper we proposed an investigation of the basic categorical and algebraic properties of quantale modules, and we showed that certain operators between objects in these categories find important applications in image processing.

We showed, in Section~\ref{nucleisec}, the properties of $\Q$-module structural closure operators and $\Q$-module transforms, and their connection with $\Q$-module morphisms. In Section~\ref{appsec} we proved that certain operators used for digital image compression and analysis are special cases of $\Q$-module transforms.

Although the results seem to be promising, especially for how easily they can be applied, we cannot pretend~--- of course~--- that the applications presented are not open to further significant developments and improvements.

In fact, as we already observed in Section~\ref{appsec}, the approach via quantale modules allowed us to group together, in a unique formal context, algorithms that act on digital images in completely different ways and have been proposed for dealing with problems different in nature. Apart from the obvious (and eternal) issue of improving the results of applications, the main open problem is the following: currently, the $\Q$-modules we really encounter in these situations are exclusively $[0,1]$-modules, a very special class of modules, hence such a formal context will be redundant from this point of view, until its applications will be extended to a wider class of tasks in data management. This is probably the most important challenge in this connection. Last, we also need to take into account that meeting this challenge would naturally give rise to a further issue, namely the necessity of numerically (or, anyhow, objectively) estimating results of the applications by introducing a sort of ``measure'' on quantale modules.

\end{document}